\documentclass[12pt]{article}
\usepackage{graphicx,amsmath,amsthm,amssymb,dsfont, mathrsfs }
\numberwithin{equation}{section}
\def\NN{\mathbb N}

\def\ZZ{\mathbb Z}
\def\PP{\mathbb P}
\def\LL{\mathbb L}

\def\FF{\mathbb F}

\def\cD{\mathcal D}

\def\cC{\mathcal C}

\def\sL{\mathscr L}

\def\cM{\mathcal M}
\def\cL{\mathcal L}
\def\fL{\mathfrak L}

\def\cL{\mathcal L}
\def\cP{\mathcal P}
\def\bcP{\boldsymbol{\cP}}

\def\cD{\mathcal D}

\def\cN{\mathcal N}

\def\fT{\mathfrak T}

\def\bT{{\bf T}}
\def\bc{{\bf c}}

\def\bx{{\bf x}}
\def\by{{\bf y}}

\def\bn{{\bf n}}

\def\bv{{\bf v}}

\def\bu{{\bf u}}
\def\bw{{\bf w}}

\def\bgamma{\boldsymbol{\gamma}}

\def\bs{{\bf 0}}
\def\b1{{\bf 1}}
\def\d1{\mathds{ 1}}
\def\dd{{\rm d}}
\def\deg{{\rm deg}}
\def\div{{\rm div}}
\def\Div{{\rm Div}}
\def\Res{{\rm Res}}
\def\Tr{{\rm Tr}}
\def\mod{{\rm mod}}
\def\ad{{\rm and}}

\def\card{{\rm card}}
\def\with{{\rm with}}
\def\where{{\rm where}}
\def\then{{\rm then}}
\def\for{{\rm for}}

\def\vol{{\rm vol}}
\def\ker{{\rm ker}}

\def\qed{ \ \vrule width.2cm height.2cm depth0cm\smallskip}
\begin{document}
%\hsize=14.8 true cm
%\vsize=30.0  true cm
%\voffset=-0.2 true cm
%\hoffset=-1.5 true cm
%%%%% To ease editing, for IMPAN journals add:

%\baselineskip=15pt

%%%%%%%%%%%%%%%%

\title{On a bounded  remainder set for $(t,s)$ sequences I}
\author{Mordechay B. Levin}

\date{}

\maketitle
\centerline{\it Dedicated to the 100th anniversary of Professor N.M. Korobov}
%$ $\\
%
\begin{abstract}
Let $\bx_0,\bx_1,...$ be a sequence of points in $[0,1)^s$.
 A subset $S$ of $[0,1)^s$ is called a bounded remainder set if there exist two real numbers $a$ and $C$ such that, for every integer $N$,
%We say that $S \subset [0,1)^s$ is a bounded remainder set with respect to the sequence $(\bx_n)_{n \geq 1}$ if there is a constant $C$ such that
$$
  | {\rm card}\{n <N \; | \;   \bx_{n} \in S\} -  a N| <C  .
$$

Let $ (\bx_n)_{n \geq 0} $ be an $s-$dimensional Halton-type sequence obtained from a global function field, $b \geq 2$,
$\bgamma =(\gamma_1,...,\gamma_s)$,
$\gamma_i \in [0, 1)$, with $b$-adic expansion $\gamma_i= \gamma_{i,1}b^{-1}+ \gamma_{i,2}b^{-2}+...$, $i=1,...,s$.
 In this paper, we prove that $[0,\gamma_1) \times ...\times [0,\gamma_s)$ is the bounded remainder set  with respect to the
sequence $(\bx_n)_{n \geq 0}$ if and only if
\begin{equation} \nonumber
   \max_{1 \leq i \leq s} \max \{ j \geq 1 \; | \; \gamma_{i,j} \neq 0 \}  < \infty.
\end{equation}
We also obtain  the similar results for a generalized Niederreiter sequences, Xing-Niederreiter sequences and  Niederreiter-Xing sequences.
\end{abstract}
%\baselineskip 15pt
%\maketitle{\bf \center }
%
%
\section{Introduction }
%${\mathscr m}, {\mathds m},{\mathcal m}, {\mathfrak m}$
%
{\bf 1.1. Bounded remainder sets.} Let $\bx_0,\bx_1,...$ be a sequence of points in $[0,1)^s$,  $S \subseteq [0,1)^s$,
\begin{equation}  \nonumber %\label{I1}
\Delta(S, (\bx_{n})_{n=0}^{N-1}  )= \sum_{n=0}^{N-1}  ( \b1_{S}(\bx_{n}) - \lambda(S)),
\end{equation}
where $\d1_{S}(\bx) =1, \; {\rm if} \;\bx  \in S$,
and $   \d1_{S}(\bx) =0,$  if $
\bx   \notin S$.  Here $\lambda(S)$ denotes the $s$-dimensional Lebesgue-measure of $S$.
We define the star {\it discrepancy} of an
$N$-point set $(\bx_{n})_{n=0}^{N-1}$ as
\begin{equation}  \nonumber % \label{I2}
   \emph{D}^{*}((\bx_{n})_{n=0}^{N-1}) =
    \sup\nolimits_{ 0<y_1, \ldots , y_s \leq 1} \; |
  \Delta([\bs,\by),(\bx_{n})_{n=0}^{N-1})/N |,
\end{equation}
where $[\bs,\by)=[0,y_1) \times \cdots \times [0,y_s) $.
 The sequence $(\bx_n)_{n \geq 0}$ is said to be
{\sf uniformly distributed}  in $[0,1)^s$
if $D_N \to \; 0$.
In 1954, Roth proved that
$ \limsup\nolimits_{N \to \infty } N (\ln N)^{-\frac{s}{ 2}} \emph{D}^{*}((\bx_{n})_{n=0}^{N-1})>0 . $
According to the well-known conjecture (see, e.g., [BeCh, p.283]), this estimate can be improved
 to
\begin{equation} \label{I2a}
 \limsup\nolimits_{N \to \infty } N (\ln N)^{-s} \emph{D}^{*}((\bx_{n})_{n=0}^{N-1})>0.
\end{equation}
See [Bi] and [Le1] for  results on this conjecture.

A sequence $(\bx_n^{(s)})_{n\geq 0}$ is of {\sf
low discrepancy} (abbreviated l.d.s.) if  $\emph{D}
((\bx_n^{(s)})_{n=0}^{N-1})\\=O(N^{-1}(\ln N)^s) $ for $ N \rightarrow \infty $.
 A sequence  of point sets $((\bx_{n,N}^{(s)})_{n=0}^{N-1})_{N=1}^{\infty}$ is of
 low discrepancy (abbreviated
l.d.p.s.) if $ \emph{D}((\bx_{n,N}^{(s)})_{n=0}^{N-1})=O(N^{-1}(\ln
N)^{s-1}) $, for $ N \rightarrow \infty $.
 For examples of such a sequences, see, e.g., [BeCh], [DiPi], and  [Ni].\\ \\
  {\bf Definition 1}. {\it Let $\bx_0,\bx_1,...$ be a sequence of points in $[0,1)^s$.
 A subset $S$ of $ [0,1)^s$ is called a {\sf bounded remainder set} for $(\bx_n)_{n \geq 0}$
if the discrepancy function $\Delta(S, (\bx_n)_{n = 0}^{N-1})$ is bounded in $\NN$.}\\

%$$  \Delta(S, (\bx_n)_{n = 1}^N) = card\{n <N \; | \;   \bx_{n} \in S -  N\lambda(S) \}$$

 Let $\alpha$ be an irrational number, let I be an interval in $[0,1)$ with  length $|I|$, let $\{n\alpha\}$ be the fractional part of $n\alpha$, $n=1,2,...$ .
 Hecke,  Ostrowski  and  Kesten  proved that $\Delta(S, (\{n\alpha\})_{n = 1}^N)$ is bounded if and only if $|I|=\{k\alpha\}$ for some integer $k$ (see references in [GrLe]).

 The sets of bounded remainder for the classical $s$-dimensional Kronecker sequence  studied
by Lev and Grepstad [GrLe]. The case of Halton's sequence was studied by Hellekalek [He].
%The references on others investigations on   bounded remainder setsee in [GrLe].

Let $b$ be a prime power,  $\bgamma =(\gamma_1,...,\gamma_s)$,
$\gamma_i \in (0, 1)$ with $b$-adic expansion $\gamma_i= \gamma_{i,1}b^{-1}+ \gamma_{i,2}b^{-2}+...$, $i=1,...,s$, and let $(\bx_n)_{n \geq 0}$ be a uniformly distributed digital Kronecker sequence.
 In  [Le1], we proved the following theorem: \\
{\bf Theorem A.} {\it The set $[0,\gamma_1) \times ...\times [0,\gamma_s)$ is of bounded remainder  with respect to  $(\bx_n)_{n \geq 0}$ if and only if}
\begin{equation} \label{Cond}
   \max_{1 \leq i \leq s} \max \{ j \geq 1 \; | \; \gamma_{i,j} \neq 0 \}  < \infty.
\end{equation}

 In this paper, we prove similar results for digital sequences described in [DiPi, Sec. 8].
 Note that according to
Larcher's conjecture [La2, p.215], the assertion of Theorem A is true for all digital $(t,s)$-sequences in base $b$.

%We will consider   in this paper some known  uniformly distributed sequences $(\bx_n)_{n \geq 1}$ satisfies the following property (condition):

\section{Definitions and auxiliary  results.}
{\bf 2.1 $(\bT, s) $ sequences.}
A subinterval $E$ of $[0,1)^s$  of the form
$$  E = \prod_{i=1}^s [a_ib^{-d_i},(a_i+1)b^{-d_i}),   $$
 with $a_i,d_i \in\ZZ, \; d_i \ge 0, \; 0  \le a_i < b^{d_i}$, for $1 \le i \le s$ is called an
{\it elementary interval in base $b \geq 2$}.\\ \\
{\bf Definition 2}. {\it Let $0 \le t \le m$  be  integers. A {\sf $(t,m,s)$-{\sf net in base $b$}} is a point set
$\bx_0,...,\bx_{b^m-1}$ in $ [0,1)^s $  such that $\# \{ n \in [0,b^m -1] | x_n \in E \}=b^t$   for every elementary interval E in base  $b$ with
$\vol(E)=b^{t-m}$.}\\  \\
{\bf Definition 3}. {\it
 Let $t \geq 0$  be an integer. A sequence $\bx_0,\bx_1,...$ of points in $[0,1)^s$ is a
{\sf $(t,s)$-{\sf sequence in base} $b$} if, for all integers $k\ge 0$ and $m \geq t$, the point set
consisting of $\bx_n$ with  $ kb^m \leq  n < (k+1)b^m$      is a $(t,m,s)$-net in base $b$.}

 By [Ni, p. 56,60], $(t,m,s)$ nets and $(t,s)$ sequences  are of low discrepancy.
See reviews
  on  $(t,m,s)$ nets and $(t,s)$ sequences  in [DiPi] and [Ni].\\  \\
%{\bf 4.2 Digital sequences and $(\bT, s) $ sequences} (\cite[Section 4]{DiPi}).\\
%
%
{\bf Definition 4.} (\cite[Definition 4.30]{DiPi}){ \it For a given dimension $s \geq 1$, an integer base $b \geq 2$, and a
function $\bT : \NN_0 \to \NN_0$ with $\bT(m) \leq  m$ for all $m \in \NN_0$, a sequence $(\bx_0,\bx_1, . . .)$
of points in $[0, 1)^s$ is called a $(\bT, s)$-sequence in base $b$ if for all integers $m \geq 0$
and $k \geq  0$, the point set consisting of the points $x_{kb^m}, . . . ,x_{kb^m+b^m-1}$ forms
a $(\bT(m),m, s)$-net in base $b$. }\\  \\
{\bf Definition 5.}   (\cite[Definition 4.47]{DiPi}){ \it
Let $m, s \geq 1$ be integers. Let $C^{(1,m)},...,C^{(s,m)}$ be $m \times m$ matrices over $\FF_b$.
Now we construct $b^m$ points in $[0, 1)^s$.
 For $ n= 0, 1,...,b^m-1$, let $n =\sum^{m-1}_{j=0} a_j(n) b^{j}$
be the $b$-adic expansion of $n$.
 Choose a  bijection
$\phi : \; \ZZ_b := \{0, 1,...., b-1\} \mapsto \FF_b$ with
$\phi(0) = \bar{0}$, the neutral element of addition in $\FF_b$.
%Let $|\phi(a)|:= |a|$ for $a \in Z_b\$.
We identify $n$ with the row vector
\begin{equation} \label{Ap300}
%Di00
                    \bn = (\bar{a}_0(n),...,\bar{a}_{m-1}(n)) \in  \FF^m_b
										\; \with \;  \bar{a}_r(n) = \phi(a_r(n)), \; r \in [0,m).
\end{equation}
We map the vectors
\begin{equation} \label{Ap301}
%\label{Di02}
	y^{(i)}_{n} =(y^{(i)}_{n,1},...,y^{(i)}_{n,m})   :=  \bn C^{(i,m)\top}, \quad
                  y^{(i)}_{n,j}=  \sum_{  r = 0}^{\infty} \bar{a}_r(n)  c^{(i)}_{j,r}     \in  \FF_b,
\end{equation}
to the real numbers
\begin{equation} \label{Ap302}
%\label{Di04}
   x^{(i)}_n =\sum_{j=1}^m  x^{(i)}_{n,j} /b^j, \quad x^{(i)}_{n,j}  =\phi^{-1} (y^{(i)}_{n,j})
\end{equation}
to obtain the point
\begin{equation} \label{Ap303}
%\label{Di06}
   \bx_n:= (x^{(1)}_n,...,x^{(s)}_n) \in [0,1)^s.
\end{equation}
 \\

%{\bf Definition 7.}  We call the point set $ \{\bx(0),...,\bx(b^m-1)\}$ constructed as introduced above a {\it digital net over} $\FF_b$ with {\it generating matrix}
%$(C^{(1)},...,C^{(s)}) $. \\

The point set  $ \{\bx_0,...,\bx_{b^m-1} \}$ is called a {\sf digital net} (over $\FF_b$) (with {\sf generating matrices} $(C^{(1,m)},...,C^{(s,m)}) $).

For $m = \infty$, we obtain a sequence $\bx_0, \bx_1,...$ of   points in $[0, 1)^s$  which is called a {\sf digital sequence} (over $\FF_b$) (with {\sf generating matrices} $(C^{(1,\infty)},...,C^{(s,\infty)}) $).}\\
We abbreviate  $C^{(i,m)}$ as $C^{(i)}$ for $m \in \NN$ and for $m=\infty$.\\ \\
{\bf 2.2  Duality theory} (see \cite[Section 7]{DiPi}).\\
Let $\cN$ be
an arbitrary $\FF_b$-linear subspace of $\FF_b^{sm}$. Let $H$ be a matrix over $\FF_b$    consisting
of $sm$ columns such that the row-space of $H$ is equal to $\cN$. Then we define
the {\sf dual space} $\cN^{\bot} \subseteq \FF_b^{sm}$
 of $\cN$  to be the null space of $H$ (see [DiPi, p. 244]). In other words,
 $\cN^{\bot} $
is the orthogonal complement of $\cN$ relative to the standard inner product
in $\FF_b^{sm}$,
\begin{equation} \nonumber % \label{Ap304}
%\label{Di10}
 \cN^{\bot} = \{ A \in  \FF_b^{ sm}  \;  | \; B \cdot A =0 \quad {\rm for \; all} \;\; B \in \cN  \}.
\end{equation}

%%%%%%%%%%%%%%%%%%%%%%%

Let  $C^{(1)},...,C^{(s)} \in \FF_b^{\infty \times \infty}$  be
generating matrices of a digital sequence  $(\bx_n(C))_{n \geq 0}$ over $\FF_b$.
For any $m \in \NN$,
we denote the $m \times m$ left-upper sub-matrix of
$C^{(i)}$ by $[C^{(i)}]_m$.
The matrices $[C^{(1)}]_m,...,[C^{(s)}]_m$
 are then the generating matrices
of a digital net. We define the {\sf overall generating matrix} of
this digital net by
\begin{equation} \label{Dif08}
  [C]_m = ([C^{(1)}]_m^{\top}|[C^{(2)}]_m^{\top}|...|[C^{(s)}]_m^{\top}) \in  \FF_b^{m \times s m}
\end{equation}
for any $m \in \NN$.

Let $\cC_m$ denote the row space of
the matrix $[C]_m$ i.e.,
\begin{equation}   \nonumber % \label{Ap309}
%\label{NiXi11}
\cC_m = \Big\{  \Big(\sum_{r=0}^{m-1}   c^{(i)}_{j,r}  \bar{a}_r(n)\Big)_{1 \leq j
\leq m, 1 \leq i \leq s}   \; | \; 0 \leq n < b^m   \Big\}.
\end{equation}
The dual space is then given by
\begin{equation} \nonumber % \label{Dif20}
 \cC_m^{\bot} = \{ A \in  \FF_b^{ sm}  \;  | \; B  \cdot A^{\top} =\bs \quad {\rm for \; all} \;\; B \in \cC_m   \}.
\end{equation}
%%%%%%%%%%%%%%%%%%%%
%
%
 {\bf Lemma A}. (\cite[Theorem 4.86]{DiPi}) {\it Let $b$ be a prime power. A strict digital
$(\bT, s)$-sequence over $\FF_b$ is
uniformly distributed modulo one, if and only if $\liminf_{m \to \infty} (m - \bT(m))=\infty$.}\\ \\
{\bf 2.3 Admissible sequences.}

 For $x =\sum_{j \geq 1}  x_{j} b^{-j}$, and $y =\sum_{j \geq 1}  y_{j}b^{-j}$
where $x_{j},y_j \in  \ZZ_b := \{0, 1,...., b\\-1\}$, we define the ($b$-adic) digital shifted point $v$ by
$v = x \oplus y := \sum_{j \geq 1}  v_{j}b^{-j}$,
 where $v_j \equiv x_j + y_j \;(\mod \;b)$ and $v_j \in \ZZ_b$. Let $\bx =
(x^{(1)},...,x^{(s)}) \in [0, 1)^s$, $\by = (y^{(1)},...,y^{(s)}) \in [0, 1)^s$.
We define the ($b$-adic) digital shifted point $\bv$ by
$ \bv =\bx \oplus \by =(x^{(1)} \oplus y^{(1)}, . . . ,x^{(s)} \oplus  y^{(s)})$.
 For $n_1,n_2 \in [0,b^m)$, we define
$n_1 \oplus n_2 :\\= (n_1 /b^m\oplus n_2)b^m)b^m$.

For $x =\sum_{j \geq 1}  x_{i}b^{-i}$,
where $x_{i} \in  \ZZ_b$,  $x_i=0 $ $(i=1,...,k)$ and $x_{k+1} \neq 0$, we define the
absolute valuation  $\left\|.  \right\|_b $ of $x$ by  $\left\|x  \right\|_b =b^{-k-1}$.
Let $\left\| n  \right\|_b =b^k$ for $n \in [b^k,b^{k+1})$.\\ \\
{\bf Definition 6.} {\it A  point set   $ (\bx_{n})_{0 \leq n <b^m} $
	in $[0,1)^s$ is  $d-${\sf admissible} in
base $b$  if}
\begin{equation} \nonumber % \label{Dif5}
  \min_{0 \leq k <n < b^m} \left\| \bx_n \ominus \bx_k  \right\|_b
  > b^{-m-d}  \quad {\rm where} \quad \left\| \bx  \right\|_b := \prod_{i=1}^s \left\|x^{(i)}_{j}  \right\|_b .
\end{equation}
{\it A sequence  $ (\bx_{n})_{n \geq 0} $
	in $[0,1)^s$ is  $d-${\sf admissible} in
base $b$ if\;} $ \inf_{n >k \geq 0}
	\left\| n \ominus k
	\right\|_b   \\ \times \left\| \bx_n \ominus \bx_k  \right\|_b  \geq b^{-d}$. \\

By [Le2], generalized Niederreiter's sequences, Xing-Niederreiter's sequences and Halton-type $(t,s)$
sequences have $d-$admissible properties.
 In  [Le2], we proved for all    $d-$admissible  digital $(t,s)$ sequences  $(\bx_n)_{ n \geq 0} $
\begin{equation}  \nonumber % \label{In10}
      \max_{1 \leq N \leq b^m}
		N \emph{D}^{*}((\bx_{n} \oplus \bw)_{0 \leq  n < N}) \geq
					 K m^{s}
\end{equation}
with some $\bw $ and $K>0$. This result supports conjecture \eqref{I2a}.\\  \\
{\bf Definition 7.} {\it A sequence   $ (\bx_{n})_{n \geq 0} $
	in $[0,1)^s$ is  {\sf weakly admissible} in
base $b$  if}
\begin{equation}   \nonumber % \label{In20}
  \varkappa_m:= \min_{0 \leq k <n < b^m} \left\| \bx_n \ominus \bx_k  \right\|_b
  >  0\quad \forall m \geq 1\; {\rm where} \;\; \left\| \bx  \right\|_b := \prod_{i=1}^s
	\left\|x^{(i)}  \right\|_b .
\end{equation}
Let $m \geq 1$, $\tau_m =[\log_q (\kappa_m)]+m $, $\bw=(w^{(1)},...,w^{(s)} )$,  $ w^{(i)} =(w^{(i)}_{1}, ..., w^{(i)}_{\tau_m})$,
\begin{equation}   \label{End}
  g_{\bw} = \{ A \geq 1 \;  | \;  x^{(i)}_{b^m A,j}=  w^{(i)}_{j}, \; j\in [1, \tau_m],\;i \in [1,s] \} \;\; \ad \;\; g_{\bw} \neq \emptyset \;\;\; \forall \; w^{(i)}_{j} \in \ZZ_b .
\end{equation} \\
{\bf Theorem B.} (see [Le3, Proposition]) {\it   Let $(\bx_n)_{n \geq 0}$ be a  uniformly distributed weakly admissible  digital $(T,s)$-sequence in base $b$, satisfying  (\ref{End}) for  all $m \geq m_0$.
Then the set $[0,\gamma_1) \times ...\times [0,\gamma_s)$ is of bounded remainder  with respect to $(\bx_n)_{n \geq 0}$ if and only if
 (\ref{Cond}) is true.}\\

{\bf 2.4 Notation and terminology for algebraic function fields.}  For
the theory of algebraic function fields, we  follow the notation and
terminology in the books  [St] and  [Sa].

 Let $b$ be an arbitrary prime power, $\FF_b$  a finite field with $b$ elements,
$\FF_b(x)$  the rational function field over $\FF_b$, and $\FF_b[x]$
  the polynomial ring over $\FF_b$.
 For $\alpha =f/g, \; f,g \in \FF_b[x]$, let
\begin{equation} \nonumber % \label{No00}
 \nu_{\infty}(\alpha) = \deg (g) - \deg (f)
\end{equation}
be the degree valuation of $\FF_b(x)$. 	
	We define the field of Laurent series as
\begin{equation} \nonumber
  \FF_b((x)) :=  \Big\{ \sum_{i = m}^{\infty} a_i x^{i}
       \; | \; m \in \ZZ, \; a_i \in \FF_b  \Big\}.
\end{equation}

 A finite extension field $F$ of $\FF_b(x)$ is called
 an algebraic function field  over $\FF_b$. Let $\FF_b$  be algebraically
closed in $F$. We express
this fact by simply saying that $F/\FF_b$  is an algebraic function field. The genus
of $F/\FF_b$ is denoted by $g$.

A place $\cP$ of $F$ is, by definition, the maximal ideal of some valuation
ring of $F$. We denote by $O_{\cP}$ the valuation ring corresponding to $\cP$ and we
denote by $\PP_F$ the set of places of $F$.
For a place $\cP$ of $F$, we write $\nu_{\cP}$ for the normalized discrete valuation of $F$
corresponding to $\cP$, and any element $t \in F  $ with $\nu_{\cP} (t) = 1$ is called a local parameter (prime element) at $\cP$.

The  field $F_{\cP}:=O_{\cP}/\cP$ is called the residue field of $F$ with respect to $\cP$.
The degree of a place
$\cP$ is defined as
$\deg(\cP) = [ F_{\cP} : \FF_b]$.
We denote by
$\Div(F)$ the set of divisors of $F/\FF_b$.

The completion of $F$ with respect to $\nu_{\cP}$ will be
denoted by $F^{(\cP)}$. Let $t$ be a local parameter of $\cP$.
Then $F^{(\cP)}$ is isomorphic to $F_{\cP} ((t))$   (see [Sa, Theorem 2.5.20]), and
 an arbitrary element $\alpha \in F^{(P)}$ can be uniquely
expanded as  (see [Sa, p. 293])
\begin{equation}   \nonumber % \label{No06}
 \alpha =  \sum_{i = \nu_{\cP}(\alpha)}^{\infty} S_i t^{i}   \quad  {\rm where}
 \quad  S_i = S_i(t,\alpha) \in F_{\cP} \subseteq F^{(P)}.
\end{equation}
The derivative $\frac{\dd \alpha}{\dd t}$, or differentiation with respect to $t$, is defined by  (see [Sa, Definition 9.3.1])
\begin{equation} \label{No08}
 \frac{\dd \alpha}{\dd t} =  \sum_{i = \nu_{\cP}(\alpha)}^{\infty} iS_i t^{i-1}   .
\end{equation}
For an algebraic function field $F/\FF_b$, we define its set of differentials (or Hasse differentials,  H-differentials) as
\begin{equation} \nonumber
                    \Delta_F = \{y \; \dd z \; | \; y \in F, \; z \;{\rm is \;a \;separating \; element  \;for} \; F/\FF_b\}
\end{equation}
(see [St, Definition 4.1.7]). \\ \\
{\bf Lemma B.} ([St, Proposition 4.1.8] or [Sa, Theorem 9.3.13]) {\it
Let $z \in F$ be separating. Then
 every differential $\gamma \in \Delta_F$
can be written uniquely as  $\gamma =
y \; \dd z$ for some $y \in  F$.}

We define the order of $\alpha \;\dd \beta$ at $\cP$ by
\begin{equation} \label{No12}
         \nu_{\cP}(\alpha \; \dd \beta) : = \nu_{\cP}(\alpha \; \dd \beta/\dd t),
\end{equation}
where $t$ is any local parameter for $\cP$  (see [Sa, Definition 9.3.8]).

 Let $\Omega_F$ be the set of all  Weil differentials of $F/\FF_b$. There exists an $F-$linear
isomorphism of the differential module $\Delta_F$  onto $\Omega_F$ (see [St, Theorem 4.3.2] or [Sa, Theorem 9.3.15]).

For $0 \neq  \omega \in  \Omega_F$, there exists a uniquely
determined divisor $\div(\omega) \in \Div(F)$.  Such a divisor $\div(\omega)$ is called a canonical
divisor of $F/\FF_b$. (see [St, Definition 1.5.11]).  For a canonical divisor $\dot{W}$, we have   (see [St, Corollary 1.5.16])
\begin{equation} \label{No14}
     \deg(\dot{W}) = 2g -2 \quad {\rm and} \quad \ell(\dot{W}) = g.
\end{equation}
 Let $ \alpha \;\dd \beta$  be a nonzero H-differential in $F$ and let $\omega$ be the corresponding
Weil differential. Then (see  [Sa, Theorem 9.3.17], [St, ref. 4.35])
\begin{equation} \label{No16}
                     \nu_{\cP} (\div(\omega) ) = \nu_{\cP}(\alpha \;\dd \beta) ,  \quad {\rm for \; all}  \quad   \cP \in \PP_F        .
\end{equation}
 Let $ \alpha \;\dd \beta$ be an H-differential, $t$  a local parameter of $\cP$, and
\begin{equation} \nonumber
     \alpha \;\dd \beta=  \sum_{i = \nu_{\cP}(\alpha)}^{\infty} S_i t^{i} \dd t \in F^{(\cP)}.
\end{equation}
Then the {\sf residue} of $ \alpha \;\dd \beta$ (see  [Sa, Definition 9.3.10) is defined by
\begin{equation} \nonumber % \label{No18}
       \Res_{\cP} ( \alpha \; \dd \beta) := \Tr_{F_{\cP}/\FF_b} (S_{-1})
			   \in \FF_b .
\end{equation}
Let
\begin{equation} \nonumber % \label{No19}
       \Res_{\cP,t} ( \alpha ) :=  \Res_{\cP} ( \alpha \dd t) .
\end{equation}
For a divisor $\cD$ of $F/\FF_b$, let $ \cL(\cD)$  denote the Riemann-Roch space
\begin{equation} \label{R-R}
           \cL(\cD) = \{ y \in F\setminus {0} \; | \; \div(y) + \cD  \geq 0 \} \cup \{0\}.
\end{equation}
Then $\cL(\cD)$ is a finite-dimensional vector space over $\FF$, and we denote its
dimension by $\ell(\cD)$. By [St, Corollary 1.4.12],
\begin{equation}  \label{No20}
                  \ell(\cD) =\{0\} \quad  \for \quad \deg(\cD) <0.
\end{equation}
{\bf Theorem C  (Riemann-Roch Theorem).}  [St, Theorem 1.5.15, and St, Theorem 1.5.17 ] {\it Let $W$ be a canonical divisor
of $F/\FF_b$. Then for each divisor $A \in \div(F)$, $\ell (A) = \deg(A) + 1 - g + \ell(W - A) $,  and }
\begin{equation} \nonumber
                     \ell (A) = \deg(A) + 1 - g  \quad {\rm for }  \quad  \deg(A) \geq  2 g -1.
\end{equation}

\section{Statements of results.}
{ \bf 3.1  Generalized Niederreiter sequence}.
In this subsection, we introduce a generalization of the Niederreiter sequence
due to Tezuka (see \cite[Section 8.1.2]{DiPi}).
By  \cite[Section 8.1] {DiPi}, the Sobol's sequence, the Faure's sequence and the original Niederreiter sequence
 are particular cases of a generalized Niederreiter sequence.

Let $b$ be a prime power and let $p_1, . . . , p_{s} \in \FF_b[x]$ be  pairwise coprime polynomials
 over $\FF_b$. Let $e_i = \deg(p_i) \geq 1$ for $1 \leq i \leq s$.
For each $j \geq  1$ and $1 \leq i \leq s$, the set of polynomials
$\{y_{i,j,k}(x) \; : \; 0 \leq  k < e_i\}$ needs to be linearly independent $(\mod \;p_i(x))$ over
$\FF_b$.
 For integers
$1 \leq i \leq s$, $j \geq 1$ and $0 \leq k < e_i$, consider the expansions
\begin{equation} \nonumber % \label{GeNi0}
  \frac{ y_{i,j,k}(x)}{
p_i(x)^{j}} = \sum_{r\geq 0} a^{(i)} (j, k, r) x^{-r-1}
\end{equation}
over the field of formal Laurent series $\FF_b((x^{-1}))$.  Then we define the matrix
 $C^{(i)} = (c^{(i)}_{ j,r})_{j \geq 1, r \geq 0}$ by
\begin{equation}  \nonumber
%\label{Ni02}
c^{(i)}_{ j,r} = a^{(i)}(Q + 1, k, r) \in \FF_b \qquad \for \qquad1 \leq i \leq s,\;  j \geq 1,\; r \geq 0,
\end{equation}
where $j -1 = Qe_i + k$ with integers $Q = Q(i, j)$ and $k = k(i, j)$ satisfying
 $0 \leq k < e_i$.

 A digital sequence $(\bx_n)_{n \geq 0}$ over $\FF_b$  generated by the matrices $C^{(1)},...,C^{(s)}$  is called a
{\sf generalized Niederreiter sequence}  (see [DiPi, p.266]). \\ \\
{\bf Theorem D.} (see [DiPi, p.266] and  \cite[Theorem 1]{Le1}) {\it The generalized Niederreiter sequence $(\bx_n)_{n \geq 0}$ with generating matrices, defined
as above, is a digital $d-$admissible (t, s)-sequence over $\FF_b$ with
$d=e_0$, $t =e_0-s$ and $e_0=e_1+...+e_s.$}

In this paper, we will consider the case where $(x,p_i)=1$ for $1 \leq i \leq s$. We will consider the general case  in [Le4].\\ \\
{\bf Theorem 1.} {\it  With the notations as above,
the set $[0,\gamma_1) \times ...\times [0,\gamma_s)$ is of bounded remainder  with respect to  $(\bx_n)_{n \geq 0}$ if and only if (\ref{Cond}) is true.} \\

{\bf 3.2 Xing-Niederreiter sequence} (see  \cite[Section 8.4 ]{DiPi}).
Let $F/\FF_b$ be an algebraic function field with full constant field $\FF_b$ and
genus $g $.
 Assume that $F/\FF_b$ has at least one
rational place $P_{\infty}$, and let $G$ be a positive divisor of $F/\FF_b$ with
$\deg(G) = 2g$ and $P_{\infty} \notin  {\rm supp}(G)$. Let $P_1, . . . , P_s$ be $s$ distinct places of $F/\FF_b$ with $P_i \neq P_{\infty}$
for $1 \leq  i  \leq  s$. Put $e_i = \deg(P_i)$ for $1 \leq i\leq  s$.

By \cite[p.279 ]{DiPi}, we have that there exists a basis   ${w_0,w_1, . . . ,w_g}$  of $\cL(G)$ over $\FF_b$ such that
\begin{equation}  \nonumber
%\label{XiNi00}
            \nu_{P_{\infty}}(w_u) = n_u \quad {\rm for} \quad 0 \leq u \leq g,
\end{equation}
where $0 = n_0 <
n_1 < .... < n_g  \leq 2g$.
For each $1 \leq i  \leq s$, we consider the chain
\begin{equation}  \nonumber
%\label{XiNi02}
          \cL(G) \subset \cL(G + P_i) \subset \cL(G + 2P_i)\subset...
\end{equation}
of vector spaces over $\FF_b$. By starting from the basis ${w_0,w_1, . . . ,w_g}$ of $\cL(G)$
and successively adding basis vectors at each step of the chain, we obtain
for each $n \in \NN$ a basis
\begin{equation}  \nonumber %\label{XiNi04}
         \{w_0,w_1, . . . ,w_g, k^{(i)}_1,k^{(i)}_2,...,k^{(i)}_{n e_i} \}
\end{equation}
of $\cL(G + n P_i)$. We note that we then have
\begin{equation}\label{XiNi06}
    k^{(i)}_j \in \cL(G + ([(j-1)/e_i+1)]P_i) \quad {\rm for} \quad1 \leq i  \leq  s \quad {\rm and } \quad j \geq 1.
\end{equation}\\
 {\bf Lemma C}. (\cite[Lemma 8.10]{DiPi}) {\it  The system $ \{ w_0, w_1,..., w_g \} \cup
 \{ k^{(i)}_j  \}_{1 \leq i \leq s, j \geq 1} $ of elements of $F$ is linearly independent over $\FF_b$. } \\ \\
%
%  ???   The following lemma was shown in [265, Lemma 2].???
%Let $z$ be a local parameter at $P_{\infty}$.
Let $z$ be  an arbitrary local parameter  at $P_{\infty}$.
For $r \in \NN_0=\NN \cup \{ 0\}$, we put
\begin{equation}\label{XiNi08}
   z_r=  \begin{cases}
        z^r \quad {\rm if} \; r \notin \{ n_0,n_1,...,n_g  \}, \\
				  w_u  \quad {\rm if} \; r= n_u \;\; {\rm for \; some} \; u \in \{0,1,...,g\}.
			\end{cases}
\end{equation}
Note that in this case $\nu_{P_{\infty}}(z_r) = r$ for all $r \in \NN_0$. For $1 \leq i \leq s$ and $j \in \NN$,
we have $k^{(i)}_j \in  \cL(G + n P_i)$ for some $n \in \NN$ and also
$P_{\infty} \notin {\rm supp}(G + nP_i)$,
hence $\nu_{P_{\infty}}(k_j^{(i)}) \geq 0$.
 Thus we have the local expansions
\begin{equation}\label{XiNi10}
 k^{(i)}_j = \sum_{  r = 0}^{\infty} a_{j,r}^{(i)} z_r  \quad {\rm for} \;\;
    1 \leq i \leq s \quad {\rm and} \;\; j \in \NN,
\end{equation}
where all coefficients $a_{j,r}^{(i)} \in \FF_b$.
Let $H_1 =\NN_0 \setminus H_2 =\{h(0),h(1),...\}$, $H_2\\ = \{n_0,n_1,...,n_g   \}$.

For $1 \leq i \leq s$ and $j \in \NN$, we now define the
sequences
\begin{equation}\label{XiNi12}
c_{j,r}^{(i)} =a_{j,h(r)}^{(i)}, \quad  \bc^{(i)}_j =(c_{j,0}^{(i)}, c_{j,1}^{(i)},...) :=(a_{j,n}^{(i)})_{n \in \NN_0 \setminus \{n_0,...,n_g \}} =(a_{j,h(r)}^{(i)})_{r \geq 0}
\end{equation}
\begin{equation}\nonumber
   = ( \widehat{a_{j,n_0}^{(i)}},a_{j,n_0+1}^{(i)},...,\widehat{a_{j,n_1}^{(i)}},
			a_{j,n_1+1}^{(i)},....,  \widehat{a_{j,n_g}^{(i)}},a_{j,n_g+1}^{(i)},....)   \in \FF_b^{\NN},
\end{equation}
where the hat indicates that the corresponding term is deleted. We define
the matrices $C^{(1)}, . . . ,C^{(s)} \in \FF^{\NN \times \NN}_b$ by
\begin{equation}\label{XiNi14}
   C^{(i)} =(\bc^{(i)}_1,\bc^{(i)}_2,\bc^{(i)}_3,...)^{\top} \quad {\rm for } \quad 1 \leq i \leq s,
\end{equation}
i.e., the vector $\bc^{(i)}_j$
 is the $j$th row vector of $C^{(i)}$ for $1 \leq i \leq s$.\\  \\
{ \bf Theorem E} (see \cite[Theorem 8.11]{DiPi} and  \cite[Theorem 1]{Le1}).  {\it With the above notations,  we have that the matrices $C^{(1)}, . . . ,C^{(s)}$
 given by (\ref{XiNi14}) are generating matrices of the Xing-Niederreiter $d-$admissible digital $(t, s)$-sequ-\\ence $ (\bx_n)_{ n \geq 0} $ with $d= e_1+...+e_s$,
$t = g +e_1+...+e_s -s $.} \\

In order to obtain the bounded remainder set property, we will take a specific local parameter $z$.
Let $P_0 \in \PP_F$, $P_0 \not\subset \{P_1,...,P_s, P_{\infty}\}$, $P_{0} \notin  {\rm supp}(G)$ and $\deg(P_0) =e_0$.
 By the Riemann-Roch  theorem, there exists a local parameter $z$ at $P_{\infty}$,  with
\begin{equation}\label{XiNi07}
             %\deg((z)_\infty)  \leq 2g+e_1 \qquad \for \qquad
             z \in \cL( (2g +1) P_0 -P_{\infty}) \setminus \cL((2g +1)P_0 -2P_{\infty}) .
\end{equation}\\
{\bf Theorem 2.} {\it  With the notations as above,
the set $[0,\gamma_1) \times ...\times [0,\gamma_s)$ is of bounded remainder  with respect to   $(\bx_n)_{n \geq 0}$ if and only if (\ref{Cond}) is true.} \\ \\
{\bf 3.3 Generalized Halton-type sequences from global function fields.}

Let $ q\geq 2 $ be an integer
 \begin{equation}\nonumber
 n=\sum_{j\geq 1}e_{q,j}(n) q^{j-1},\;  e_{q,j}(n) \in \{0,1, \ldots
 ,q-1\}, \;  {\rm and}    \; \varphi_q(n)= \sum_{j\geq 1}e_{q,j}(n) q^{-j}.
 \end{equation}
Van der Corput    proved that $ (\varphi_q(n))_{n\geq 0}$ is a $1-$dimensional l.d.s.  Let
\begin{equation}\nonumber
  \hat{H}_s(n)= (\varphi_{\hat{q}_1}(n),\ldots ,\varphi_{\hat{q}_s}(n)), \quad n=0,1,2,...,
\end{equation}
 where $
\hat{q}_1,\ldots ,\hat{q}_s\geq 2 $ are pairwise coprime integers.
Halton   proved that $ ( \hat{H}_s(n))_{n\geq 0}$ is an $s-$dimensional l.d.s. (see [Ni]).

Let $Q=(q_1,q_2,....)$ and $Q_j=q_1q_2....q_j$, where $q_j\ge2$ \ $(j=1,2,\dots)$ is a sequence of integers.
Every nonnegative integer $n$ then has a unique $Q$-adic representation of the form
\begin{equation} \nonumber
    n= \sum_{j=1}^{\infty} n_j q_1 \cdots q_{j-1} = n_1 +n_2 q_1 + n_3 q_1 q_2 + \cdots,
\end{equation}
where $n_j \in \{ 0,1,...,q_j-1 \}$. We call this the Cantor expansion of $n$ with respect to the base $Q$.
Consider Cantor's expansion of $x \in [0,1):$
$$x=\sum\nolimits_{j=1}^{\infty} x_j / Q_j, \quad
x_j \in \{0,1,\dots,q_j-1\}, \quad x_j \neq q_j-1\; {\rm for \;infinitely \; many} \;j.$$
The $Q-$adic representation of $x$ is then unique.
We define the radical inverse function
\begin{equation} \nonumber
  \varphi_Q\Big( \sum_{j=1}^{\infty} n_j q_1 \cdots q_{j-1} \Big) = \sum_{j=1}^{\infty} \frac{n_j}{ q_1 \cdots q_{j}} .
\end{equation}
Let   $p_{i,j} \geq 2$ be integers $( s \geq i \geq 1, j \geq 1)$,  $g.c.d.(p_{i,k},p_{j,l})=1$ for $i\neq j$, $\tilde{P}_{i,0} =1, \;\;  \tilde{P}_{i,j} =\prod_{ 1 \leq k \leq j}  p_{i,k},
    \; i \in [1,s],\; j \geq 1$,
 $\cP_i = (p_{i,1},p_{i,2},...)$,  $\bcP=(\cP_1,...,\cP_s)$.
In [He], Hellecaleq proposed the following generalisation of the Halton sequence:
\begin{equation} \label{He}
    H_{\bcP}= (\varphi_{\cP_1}(n),\ldots ,\varphi_{\cP_s}(n))_{n=0}^{\infty}.
\end{equation}

In [Te], Tezuka introduced a polynomial arithmetic analogue of the Halton sequence :

Let $p(x)$ be an arbitrary nonconstant polynomial over $\FF_b$, $e=\deg(p)$, $n= a_0(n) +a_1(n)b + \cdots +a_{m}(n)b^m$. We fix
a bijection $ \phi \;: \; \ZZ_b \to \FF_b$ with $\phi(0) = \bar{0}$.
Denote $v_n(x) = \bar{a}_0(n) +\bar{a}_1(n) x + \cdots + \bar{a}_m(n) x^{m}$, where $\bar{a}_r(n) = \phi(a_r(n))$, $r=0,1,...,m$.
Then $v_n(x)$ can be represented in terms of $p(x)$ in the following way:
\begin{equation} \nonumber
v_n(x) = r_0(x) + r_1(x) p(x) + \cdots  + r_{k} (p(x))^{k}, \quad \with \quad k=[m/e].
 \end{equation}
We define the radical inverse function $\varphi_{p(x)}: \FF_b[x] \to \FF_b(x)$ as follows
\begin{equation} \nonumber
  \varphi_{p(x)}(v_n(x)) = r_0(x)/p(x) + r_1(x)/ p^2(x) + \cdots  + r_{k}/ (p(x))^{k+1}.
 \end{equation}
Let $p_1(x),..., p_s(x)$ be pairwise coprime. Then Tezuka's sequence is defined as follows
\begin{equation} \nonumber
    \bx_n= (\sigma_1(\varphi_{p_1(x)}(n)),\ldots ,\sigma_s(\varphi_{p_s(x)} (n)))_{n=0}^{\infty},
\end{equation}
where each $\sigma_i$ is a mapping from $F$ to the real field defined by
 $ \sigma_i(\sum_{j \geq w} \dot{a}_j x^{-j})\\ =  \sum_{j \geq w}  \phi^{-1}(\dot{a}_j) b^{-j}$.
By [Te], $(\bx_n)_{n \geq 0}$ is a $(t,s)$ sequence in base $b$.\\

In 2010, Levin [Le1] and in 2013, Niederreiter and Yeo  [NiYe] generalized Tezuka's construction to the case of arbitrary algebraic function fields $F$.
 The construction of [NiYe] is follows:

Let $F/\FF_b$ be an algebraic function field with full constant field $\FF_b$ and
genus $g$.
We assume that $F/\FF_b$ has at least one rational place, that is, a place of degree 1. Given a dimension $s \geq 1$, we choose $s +1$ distinct places  $P_1$,...,$P_{s}$, $P_{\infty}$ of $F$ with deg$(P_{\infty}) = 1$. The
degrees of the places $P_1$,...,$P_s$ are arbitrary and we put $e_i = \deg (P_i)$ for
 $1 \leq i \leq s$. Denote
by $O_F$ the holomorphy ring given by
$ O_F = \bigcap_{ P  \neq P_{\infty}} O_P $,
where the intersection is extended over all places $ P  \neq P_{\infty} $ of $F$, and $O_P$ is the valuation ring of $P$.
We arrange the elements of $O_F$ into a sequence by using the fact that $O_F = \bigcup_{m \geq 0} \cL(mP_{\infty}) $.
The terms of this sequence are denoted by $f_0, f_1, . . .$ and they are obtained as follows. Consider the chain  $ \cL(0) \subseteq L(P_{\infty}) \subseteq L(2P_{\infty}) \subseteq \cdots	$
of vector spaces over $\FF_b$. At each step of this chain, the dimension either remains the same
or increases by $1$. From a certain point on, the dimension always increases by $1$ according to
the Riemann-Roch theorem. Thus we can construct a sequence $v_0, v_1, . . .$ of elements of $O_F$
such that
$  \{v_0, v_1, . . . , v_{\ell(mP_{s+1})-1} \} $
is a $\FF_b$-basis of $\cL(m P_{s+1})$.  We fix
a bijection $ \phi \;: \; \ZZ_b \to \FF_b$ with $\phi(0) = \bar{0}$. Then we define
\begin{equation}\nonumber
 f_n = \sum_{  r = 0}^{\infty} \bar{a}_r(n) v_r  \in O_F\quad {\rm with} \quad
             \bar{a}_r(n)=\phi(a_r(n))
\quad {\rm for} \;\;  n=\sum_{  r =0}^{\infty} a_r(n) b^r \; .
\end{equation}
Note that the sum above is finite since for each $n \in \NN$. We have $a_r(n) = 0$ for all sufficiently
large $r$.
 By the Riemann-Roch theorem, we have
\begin{equation}  \nonumber
 \{\tilde{f} \; |\; \tilde{f} \in \cL((m+g-1) P_{s+1}) \} =
 \{f_n \; |\; n \in [0,b^m) \}
\quad \for \quad m \geq g.
\end{equation}
 For each $i = 1, . . . , s$, let $\wp_i$ be the maximal ideal of $O_F$ corresponding to $P_i$. Then the
residue class field $F_{P_i} := O_F /\wp_i$ has order $b^{e_i}$ (see [St, Proposition 3.2.9]). We fix a bijection
$   \sigma_{P_i}  \; : \;F_{P_i} \to Z_{b^{e_i}} $.
 For each
$i = 1, . . . , s$, we can obtain a local parameter $t_i \in O_F$ at $\wp_i$, by applying the Riemann-Roch
theorem and choosing
$ t_i \in  \cL(k P_{\infty} - P_i) \setminus  \cL(k P_{\infty} - 2P_i) $
for a suitably large integer $k$. We have a local expansion of
$f_n$ at $\wp_i$ of the form
\begin{equation}\nonumber
  f_n = \sum_{  j \geq 0}  f^{(i)}_{n,j} t_i^j    \quad {\rm with \; all} \;\;   f^{(i)}_{n,j}    \in F_{P_i} ,\; n=0,1,...\; .
\end{equation}
%Compare with [Di, pp. 604-605] for a detailed description of how to obtain a local expansion.
We define the map $ \xi  \; : \; O_F  \to [0, 1]^s$ by
\begin{equation}\nonumber
       \xi(f_n) = \Big(  \sum_{  j = 0}^{\infty}  \sigma_{P_1} ( f^{(1)}_{n,j}) (b^{e_1})^{-j-1},...,
						   \sum_{  j = 0}^{\infty}  \sigma_{P_s} ( f^{(s)}_{n,j}) (b^{e_s})^{-j-1} \Big) .	
\end{equation}
Now we  define the sequence $\bx_0, \bx_1, . . .$ of points in $[0, 1]^s$ by
$  \bx_n  = \xi(f_n)   \;\;\;\; {\rm for} $\\$n=0,1,...$ .
 From [NiYe, Theorem 1], we get the following theorem :\\ \\
{\bf Theorem F.}  {\it With the notation as above, we have that $(\bx_n)_{n \geq 0} $
%is a $(t,m,s)$-net
is a  $(t, s)$-sequence over $\FF_b$ with
 $t=g+ e_1+...+e_s-s$.}\\

The construction of Levin [Le1] is similar, but more complicated than in [NiYe].
However in [Le1], we can use arbitrary pairwise coprime divisors $D_1,...,D_s$  instead of places $P_1,...,P_s$.\\

In this paper, we introduce the Hellecalek-like generalisation \eqref{He} of the above construction:

Let $\PP_F :=\{P | P \; {\rm be\; a \;place \;of\;} F/\FF_b\}$,  $P_0,P_{\infty} \in \PP_F$, $\deg(P_{\infty})=1$, $\deg(P_{0})\\=e_0$, $P_0 \neq P_{\infty} $,
 $P_{i,j} \in \PP_F$ for $1 \leq j , 1 \leq i \leq s$, $P_{i_1,j_1} \neq P_{i_2,j_2}$ for $i_1 \neq i_2$,
$P_{i,j} \neq P_0$, $P_{i,j} \neq P_{\infty}$ for all $i,j$, $\dot{n}_{i,j}= \deg(P_{i,j})$,
$ n_{i,j}=\deg(\cP_{i,j}), \; \cP_{0,j} =P_0^j$,
\begin{equation} \label{T3-0}
 \cP_{i,0} =1, \;  \cP_{i,j} =\prod_{ 1 \leq k \leq j}  P_{i,k}, \; n_{i,j} =\deg(\cP_{i,j})
  = n_{i,j-1} + \dot{n}_{i,j}, n_{i,0}=0.
\end{equation}

Let $i \in [0,s]$. We will construct a basis $(w^{(i)}_{j})_{j \geq 0}$ of $O_F$ in the following way. Let
\begin{align}
&   L_{i,j}=\cL((n_{i,j} +2g-1)P_{\infty}) =\cL(A_{i,j}), \;  A_{i,j} = (n_{i,j} +2g-1)P_{\infty},
              \label{T3-1}  \\
&   \fL_{i,j}=\cL((n_{i,j} +2g-1)P_{\infty}- \cP_{i,j}) = \cL(B_{i,j}), \;    B_{i,j} = (n_{i,j} +2g-1)P_{\infty}- \cP_{i,j},   \nonumber   \\
&   \sL_{i,j}=\cL((n_{i,j} +2g-1)P_{\infty}- \cP_{i,j-1})     , \;\;   \dot{B}_{i,j} =(n_{i,j} +2g-1)P_{\infty}- \cP_{i,j-1}. \nonumber
\end{align}
Using the Riemann-Roch theorem,  we obtain
\begin{align}
  \deg(A_{i,j})= n_{i,j} +2g-1, \;\; \dim(L_{i,j})= n_{i,j} +g, \;\;
  \deg(B_{i,j})= 2g-1,   \label{T3-2}     \\
 \dim(\fL_{i,j})= g, \quad  \deg(\dot{B}_{i,j})= \dot{n}_{i,j} +2g-1, \quad \dim(\sL_{i,j})= \dot{n}_{i,j} +g.  \nonumber
\end{align}
Let $(u^{(i)}_{j,\mu})_{\mu =1}^{g}$ be a $\FF_b$ linear basis of $\fL_{i,j}$.
By \eqref{T3-1} and \eqref{T3-2}, we get that the basis  $(u^{(i)}_{j,\mu})_{\mu =1}^{g}$ can be extended to a basis
$(v^{(i)}_{j,1}, \cdots, v^{(i)}_{j, \dot{n}_{i,j} }, u^{(i)}_{j,1}, \cdots , u^{(i)}_{j,g})$
 of $\sL_{i,j}$.

  Bearing in mind that  $(u^{(i)}_{j,\mu})_{\mu =1}^{g}$ is a $\FF_b$ linear basis of $\fL_{i,j}$, we obtain that $v^{(i)}_{j,\mu} \notin \fL_{i,j}$ for $ \mu \in [1, \dot{n}_{i,j}  ]$.
So
\begin{equation} \label{T3-3}
     v^{(i)}_{j,\mu} \in \LL_{i,j}:= \sL_{i,j} \setminus \fL_{i,j} \quad \for \quad \mu \in [1, \dot{n}_{i,j}  ].
\end{equation}
Let
\begin{equation} \label{T3-3a}
   V_{i,j}:= \{ v^{(i)}_{k,\mu}  \; | \; 1 \leq \mu \leq \dot{n}_{i,k},\; 1 \leq k \leq j \}  \cup \{ u^{(i)}_{j,\mu} \; | \; \mu =1,...,g \}.
\end{equation}
We claim that vectors from $V_{i,j}$ are $\FF_b$ linear independent. Suppose the opposite. Assume that there exists $ b^{(i)}_{k,\mu} \in \FF_b $ such that
\begin{equation} \label{T3-4}
    \dot{u} + \ddot{u} =0, \;\; \where \;\;   \dot{u} = \sum_{k=1}^j w_k, \;\;
    w_k  =\sum_{\mu=1}^{\dot{n}_{i,k}}  b^{(i)}_{k,\mu} v^{(i)}_{k,\mu}, \;\;
     \ddot{u} = \sum_{\mu=1}^{g}  b^{(i)}_{0,\mu} u^{(i)}_{j,\mu} .
\end{equation}
%and $ \ddot{u} = \sum_{\mu=1}^{g}  b^{(i)}_{0,\mu} u^{(i)}_{j,\mu}   $.
%
Let $w_{l} \neq 0$ for some $l \in [1,j]$  and let $w_k=0$ for all $k \in [1,l) $. \\
Using \eqref{T3-1} - \eqref{T3-3}, we get
\begin{equation} \nonumber %   \label{T3-5}
 w_l \in \LL_{i,l}=\cL((n_{i,l} +2g-1)P_{\infty}- \cP_{i,l-1}) \setminus  \cL((n_{i,l} +2g-1)P_{\infty}- \cP_{i,l}) .
\end{equation}
Applying definition \eqref{R-R} of  the Riemann-Roch space, we obtain
\begin{equation} \nonumber %\label{T3-6}
 w_l \in \cL((n_{i,j} +2g-1)P_{\infty}- \cP_{i,l-1}) \setminus  \cL((n_{i,j} +2g-1)P_{\infty}-\cP_{i,l}) .
\end{equation}
But from  \eqref{T3-4}, \eqref{T3-0} and \eqref{T3-3}, we have
\begin{equation} \nonumber
   -w_l = \dot{u} + \ddot{u} -\sum_{k=1}^l w_k =   \sum_{k=l+1}^j w_k  + \ddot{u}
   \in \cL((n_{i,j} +2g-1)P_{\infty}- \cP_{i,l}) .
\end{equation}
We have a contradiction. Hence  vectors from $V_{i,j}$ are $\FF_b$ linear independent.  \\
By   \eqref{T3-1}  - \eqref{T3-3a}, we have $V_{i,j} \subset L_{i,j}$ and
\begin{equation} \nonumber           %\label{T3-8}
 \card(V_{i,j}) =\sum_{k=1}^j {\dot{n}_{i,k}} +g = n_{i,j}+g = \dim(L_{i,j}).
\end{equation}
Hence vectors from $V_{i,j}$ are the $\FF_b$ linear basis of $ L_{i,j} $.

Now we will find a basis of $L_{i, j-2g}$. We claim that $  u^{(i)}_{j,\mu} \notin   L_{i, j-2g}$
 for $\mu \in [1,g]$. Suppose the opposite. By    \eqref{T3-1} and \eqref{T3-2}, we get
\begin{align}    %\label{T3-9}
 &  u^{(i)}_{j,\mu} \in   L_{i, j-2g} \cap \fL_{i, j} =   \cL((n_{i,j-2g} +2g-1)P_{\infty} )
   \cap \cL((n_{i,j} +2g-1)P_{\infty}- \cP_{i,j})  \nonumber \\
 &   = \cL((n_{i,j-2g} +2g-1)P_{\infty}- \cP_{i,j}) =
   \cL(T) .  \nonumber
\end{align}
By    \eqref{T3-0}, $\deg(T) =n_{i,j-2g} +2g-1 -n_{i,j} <0 $. Hence $ \cL(T) = \{ 0 \}$. We have a contradiction.
Bearing in mind that $V_{i,j}$  is $\FF_b$ linear basis of $L_{i,j} $, we obtain that a basis of $L_{i, j-2g}$ can be chosen from the set $(v^{(i)}_{1,1},...,v^{(i)}_{1,\dot{n}_{i,1}}, ...,   v^{(i)}_{j,1},..., v^{(i)}_{j,\dot{n}_{i,j}})$\\ $= V_{i,j}\setminus  \{ u^{(i)}_{j,\mu} \; | \; \mu =1,...,g \} $.
From   \eqref{T3-1}  - \eqref{T3-3}, we get
\begin{equation} \nonumber %   \label{T3-10}
     v^{(i)}_{k,\mu} \in \sL_{i,k} \subseteq L_{i,j-2g}  \quad \for \quad \mu \in [1, \dot{n}_{i,k}] \;\; \ad \;\; 1 \leq k \leq j-2g.
\end{equation}
Hence vectors
\begin{equation} \nonumber %\label{T3-11}
 v^{(i)}_{1,1},...,v^{(i)}_{1,\dot{n}_{i,1}}, ...,   v^{(i)}_{j-2g,1},..., v^{(i)}_{j-2g,\dot{n}_{i,j-2g}},
\tilde{v}^{(i)}_{j,1},...,\tilde{v}^{(i)}_{j,g}
      \quad \with \quad    \tilde{v}^{(i)}_{j,\mu} =v^{(i)}_{k,\rho},
\end{equation}
    $ 1 \leq \mu \leq g$,  for some  $ \rho \in [1, \dot{n}_{i,k} ]$  and $  k \in (j-2g, j]$  are an $\FF_b$ linear basis of $L_{i,j -2g}$  $(0 \leq i \leq s)$.

Therefore $(v^{(i)}_{k,\mu})_{1 \leq \mu \leq \dot{n}_{i,k, k \geq 1}} $ is  the $\FF_b$ linear basis of $O_F= \cup_{j \geq 1} L_{i,j}$.
We put in order the basis $(v^{(i)}_{k,\mu})_{1 \leq \mu \leq \dot{n}_{i,k}, k \geq 1} $   as follows
\begin{equation}\label{T3-12}
   w^{(i)}_{n_{i,j-1} + \mu -1} = v^{(i)}_{j,\mu}, \quad \with \quad n_{i,0}=0,\;
     1 \leq \mu \leq  \dot{n}_{i,j}, \; 0 \leq i \leq s.
\end{equation}
So we proved the following lemma : \\  \\
{\bf Lemma 1.} { \it For all $i \in [0,s]$  there exists a sequence  $(w^{(i)}_j)_{j \geq 0}$ such that
 $(w^{(i)}_{j})_{j \geq 0}$  is  a $\FF_b$ linear basis of $O_F$ and for all $j \geq 1$
 a $\FF_b$ linear basis of $L_{i,j}$ can be chosen from the set $\{ w^{(i)}_{0},..., w^{(i)}_{n_{i,j+2g}-1}\}$.}  \\

Bearing in mind that $(w^{(i)}_{j})_{j \geq 0}$  is the $\FF_b$ linear basis of  $O_F$, we obtain for all $i \in [1,s]$ and $r \geq 0$
 that there exists $c^{(i)}_{j,r} \in \FF_b$ and integers $l^{(i)}_r$
 such that
\begin{equation} \label{T3-13}
%\label{Cond1}
   w^{(0)}_{r} = \sum_{j=1}^{l^{(i)}_r }   c^{(i)}_{j,r} w^{(i)}_{j-1}, \quad c^{(0)}_{j,j-1}=1,\; \ad \;
                c^{(0)}_{j,r} =0 \; \for \; j-1\neq r.
\end{equation}
Let $  n = \sum_{  r \geq 0} a_r(n) b^r$. We fix
a bijection $ \phi \;: \; \ZZ_b \to \FF_b$ with $\phi(0) = \bar{0}$. Then we define
\begin{equation}        \label{T3Pr0a}
 f_n = \sum_{  r = 0}^{\infty} \bar{a}_r(n) w^{(0)}_{r } \in O_F \;\;  {\rm with} \;\;
             \bar{a}_r(n)=\phi(a_r(n))
\;\; {\rm for} \;\;  n=0,1,... \; .
\end{equation}
By \eqref{T3-13}, we have for $i \in [0, s]$
\begin{equation}\label{T3Pr1}
 f_n = \sum_{  r = 0}^{\infty} \bar{a}_r(n) \sum_{  j = 1}^{l_{i,r}} c^{(i)}_{j,r} w^{(i)}_{j-1} =
           \sum_{  j =1}^{\infty} w^{(i)}_{j-1}  \sum_{  r = 0}^{\infty} \bar{a}_r(n)  c^{(i)}_{j,r}=
         \sum_{  j =1}^{\infty}  y^{(i)}_{n,j}  w^{(i)}_{j-1}
\end{equation}
where $ y^{(i)}_{n,j}=  \sum_{  r \geq 0} \bar{a}_r(n)  c^{(i)}_{j,r}     \in  \FF_b$,  $y^{(0)}_{n,j}=\bar{a}_{j-1}(n) $.\\
We map the vectors
\begin{equation} \label{T3Pr2}
%\label{Di02}
	y^{(i)}_{n} =(y^{(i)}_{n,1},y^{(i)}_{n,2},...)
\end{equation}
to the real numbers
\begin{equation} \nonumber %   \label{T3Pr3}
%\label{Di04}
   x^{(i)}_n =\sum_{j \geq 1} \phi^{-1} (y^{(i)}_{n,j})/b^j% \quad
%   x^{(0)}_n =\sum_{j \geq 1}   a_{j-1}(n)/b^j
\end{equation}
to obtain the point
\begin{equation} \label{T3Pr4}
%\label{Di06}
   \bx_n:= (x^{(1)}_n,...,x^{(s)}_n) \in [0,1)^s.
\end{equation} \\
{\bf Theorem 3.} {\it  With the notations as above,
the set $[0,\gamma_1) \times ...\times [0,\gamma_s)$ is of bounded remainder  with respect to $(\bx_n)_{n \geq 0}$ if and only if (\ref{Cond}) is true.} \\  \\
{\bf Remark.} It is easy to verify that Hellekalek's sequence and our generalized Halton-type sequence  $(\bx_n)_{n \geq 0}$ are  l.d.s if
\begin{equation}  \nonumber
  \limsup_{m \to \infty} m^{-s}\sum_{i=1}^s \sum_{j=1}^m \log(p_{i,j}) < \infty  \;\; \ad \;\;
    \limsup_{m \to \infty} m^{-s}\sum_{i=1}^s  \sum_{j=1}^m  \deg(P_{i,j}) < \infty.
\end{equation} \\

%%%%%%%%%%%%%%%%%%%%%%%%%%%%%%%%%%%%%%%%%%%%%%%%%%%%%%%%%%%%%%%%%%%%%%%%%%%%%%%%%%%%%%%%%%%%%%%%%%%%%%%%%%%%%%%%%%%%%%%%%%%%%%%%%%%%%%%%%%%%%%%%%%%%%%%%%%%%%%%%%%%%%%%%%%%%%%%%%%%%%%%%%%%%

{\bf 3.4 Niederreiter-Xing sequence}  (see  \cite[Section 8.3 ]{DiPi}).
Let $F/\FF_b$ be an algebraic function field with full constant field $\FF_b$ and
genus $g$.
 Assume that $F/\FF_b$ has at least $s+1$ rational places. Let $ P_1,...,P_{s+1}$ be $s+1$
distinct rational places of $F$. Let  $G_m =m(P_1+...+P_s) -(m-g+1)P_{s+1}$, and
let $t_i$ be a local parameter at $P_i$, $1 \leq i \leq s+1$.
For any $f \in \cL(G_m)$ we have $\nu_{P_i}(f) \geq  - m$, and so the local
expansion of $f$ at $P_i$ has the form
\begin{equation}  \nonumber
%\label{NiXi00}
f = \sum_{j=-m}^{\infty} f_{i,j} t^j_i , \quad {\rm with } \quad
  f_{i,j}  \in  \FF_b, \; j \geq - m, \; 1 \leq i \leq s.
\end{equation}
For $1 \leq i \leq s$, we define the $\FF_b$-linear map $\psi_{m,i} \;:\; \cL(G_m) \to \FF^m_b$
by
\begin{equation} \nonumber
%\label{NiXi01}
             \psi_{m,i}(f)= (   f_{i,- 1}, . . . , f_{i,- m}) \in \FF^m_b,\quad \for \quad
     f \in \cL(G_m).
\end{equation}
Let
\begin{equation}  \nonumber
%\label{NiXi03}
        \cM_m=   \cM_m(P_1,...,P_s;G_m):= \{   (\psi_{m,1}(f), . . . , \psi_{m,s}(f))  \in \FF^{ms}_b\;  |  \;f \in \cL(G_m) \}.
\end{equation}

Let  $C^{(1)},...,C^{(s)} \in \FF_b^{\infty \times \infty}$  be the
generating matrices of a digital sequence  $\bx_n(C)_{n \geq 0}$, and let $(\cC_m )_{m \geq 1}$ be the associated sequence of
 row spaces of overall
generating matrices $[C]_m$, $m=1,2,...$ (see (\ref{Dif08})).
  \\ \\
{ \bf Theorem G.} (see \cite[Theorem 7.26 and Theorem 8.9]{DiPi})  {\it There exist matrices $C^{(1)},...,C^{(s)}$
such that  $(\bx_n(C))_{n \geq 0}$ is a digital $(t, s)$-sequence  with $t = g$ and $\cC_m^{\bot} =\cM_m(P_1,...,P_s; G_m)$ for $m \geq g+1$, $s \geq 2$.} \\

In [Le2, p.24], we proposed the following way to get $\bx_n(C)_{n \geq 0}$ : \\

We consider the $H$-differential  $dt_{s+1}$. Let $\omega$ be the corresponding
Weil differential,  $\div(\omega)$ the divisor of $\omega$, and $W:=\div(dt_{s+1})= \div(\omega)$.
By  (\ref{No08})-(\ref{No14}), we have $ \deg(W)=2g-2$.
We consider  a sequence $\dot{v}_0, \dot{v}_1, . . .$ of elements of $F$ such that $ \{\dot{v}_0, \dot{v}_1, . . . , \dot{v}_{\ell((m -g+1)P_{s+1} +W)-1} \}$
is an $\FF_b$ linear basis of $L_m:=\cL((m -g+1)P_{s+1} +W)$ and
\begin{equation}\label{NiXi11a}
 \dot{v}_{r} \in L_{r+1} \setminus L_{r}, \; \quad  \nu_{P_{s+1}}(  \dot{v}_{r}) =-r+g-2,\;
  r \geq g,\;\; \ad \;\; \dot{v}_{r,r+2-g} =1, \; \dot{v}_{r,j} =0
\end{equation}
for $2  \leq j <r +2- g$, where
\begin{equation} \nonumber
  \dot{v}_r   :=  \sum_{j \leq r-g+2}
 \dot{v}_{r,j}t_{s+1}^{-j}   \quad {\rm for}  \quad
	\dot{v}_{r,j} \in \FF_b \;\; \ad \;\; r \geq g.
\end{equation}
According to Lemma B, we have that there exists $\tau_i \in F$
$(1 \leq i \leq s)$ such that $ \dd t_{s+1} = \tau_i \dd t_i,   \quad \for  \quad 1 \leq i \leq s$.\\
%\begin{equation}\label{NiOz23d}
%   \dd t_s = \tau_i \dd t_i, \qquad 1 \leq i \leq s.
%\end{equation}
%Using (\ref{Hal03}), we define $\dot{v}_j =v_j \tau_0$, for $j \geq 0$. According to (\ref{NiXi04}), we obtain
%$v_j \in \cL((m+g-1) P_{s+1} )$ if and only if $\dot{v}_j \in \cL((m-g+1) P_{s+1} +W)$.
Bearing in mind (\ref{No12}), (\ref{No16}) and (\ref{NiXi11a}),  we get
%Using (\ref{No12}) and (\ref{NiOz23c}),   we get
\begin{equation} \nonumber
% \label{NiXi13}
   \nu_{P_i}( \dot{v}_j \tau_i) =\nu_{P_i}( \dot{v}_j \tau_i \dd t_i) =\nu_{P_i}( \dot{v}_j
	\dd t_{s+1}) \geq
		\nu_{P_i}( \div(\dd t_{s+1}) - W ) =0, \quad  j \geq 0.
\end{equation}
We consider the following local expansions
\begin{equation} \label{NiXi05}
  \dot{v}_r   \tau_i\; := \; \sum_{j = 1}^{\infty}
	\dot{c}^{(i)}_{j,r} t_i^{j-1}  , \quad {\rm where \; all }  \quad
	\dot{c}^{(i)}_{j,r} \in \FF_b, \; 1 \leq i \leq s, \; j \geq 1.
\end{equation}
%
%Applying (\ref{NiXi05}), we have \begin{equation}\label{NiXi12}    \dot{v}_j   \tau_i\; := \; \sum_{j = \dot{a}(i,j,\tau) }^{\infty}	\ddot{c}^{(i)}_{r,j} t_i^{j}  , \quad {\rm where \; all }  \quad   \ddot{c}^{(i)}_{r,j}	\in \FF_b, \; 1 \leq i \leq s, \; j \geq 0.\end{equation}
%
Now let $\dot{C}^{(i)} =(\dot{c}^{(i)}_{j,r}  )_{j-1,r \geq 0}$, $1 \leq i \leq s$,
 and let $(\dot{\cC}_m^{\bot} )_{m \geq 1}$ be the associated sequence of
 row spaces of overall
generating matrices $[\dot{C}]_m$, $m=1,2,...$ (see (\ref{Dif08})). \\  \\
{\bf Theorem H} (see [Le2, Theorem 5]). {\it  With the above notations,
$(\bx_n(\dot{C}))_{n \geq 0}$ is a digital $d-$admissible $(t,s)$ sequence
  with $d=g+s$,
$t=g$, and
 $\dot{\cC}_m^{\bot} =\cM_m(P_1,...,P_s; G_m)$ for all $m \geq g+1$.} \\

We note that condition \eqref{NiXi11a} is required in the proof of Theorem H only in order to get the discrepancy lower bound. While the equality  $\dot{\cC}_m^{\bot} =\cM_m(P_1,...,P_s; G_m)$ is true for arbitrary sequence  $\dot{v}_0, \dot{v}_1, . . .$ of elements of $\FF_b$ such that
for all $m \geq 1$
\begin{equation}\label{Th4-50}
   \{\dot{v}_0, \dot{v}_1, . . . , \dot{v}_{\ell((m -g+1)P_{s+1} +W)-1} \} \;\;
{\rm is\; a}  \;\; \FF_b \;\; {\rm linear\;  basis \;of}\; L_m.
\end{equation}

In order to obtain the bounded remainder property, in this paper, we will construct from $(\dot{v}_n)_{n \geq 0}$ a special basis $(\ddot{v}_n)_{n \geq 0}$ as follows : \\

Let $P_0 \in \PP_F$, $P_0 \neq  P_i$ $(i=1,...,s+1)$, and let $t_0 $ be a local parameter of $P_0$. For simplicity, we suppose that $\deg(P_0) =1$.
Let
\begin{align}  \nonumber  %   \label{NiXi10}
 &L_{m} =\cL((m- g+1)P_{s+1}+W), \quad     \sL_{m} =  \cL((m +2)P_{s+1}+W- m P_{0}),\\
 & \fL_{m} =\cL((m +2)P_{s+1}+W- (m+1)P_{0}).  \label{NiXi10a}
\end{align}
It is easy to verify that
\begin{align}  \nonumber  %
& \deg( \sL_{m}) =2g, \;\;
\dim( \sL_{m}) =g+1,\;\;\deg( \fL_{m}) =2g-1,\;\; \dim( \fL_{m}) =g, \\
&   \for \; m  \geq 0,\; \deg( L_{m}) =m+g-1,\;\; \dim( L_{m}) =m,\;
\for \; m \geq g \; . \label{NiXi10}
\end{align}
Using the Riemann-Roch theorem,  we have that there exists
\begin{equation}\label{NiXi11}
 w_{m} \in \sL_{m} \setminus \fL_{m}, \; \quad \ad \quad \;\;  w_{m} \in L_{m+g+1},  \;\; \quad  m=0,1,... \; .
\end{equation}
According to Lemma B, we have that there exists $\tau_0 \in F$,
  such that $ \dd t_{s+1} = \tau_0 \dd t_0$.

Let $u \in L_m= \cL((m -g+1)P_{s+1} +W)$  with $m \geq 0$.
Bearing in mind (\ref{No12}), (\ref{No16}), (\ref{NiXi10a})-(\ref{NiXi11}) and the Riemann-Roch theorem,  we get
%Using (\ref{No12}) and (\ref{NiOz23c}),   we get
\begin{equation}  \label{NiXi13}
   \nu_{P_0}( u \tau_0) =\nu_{P_0}( u \tau_0 \dd t_0) =\nu_{P_0}( u
	\dd t_{s+1}) = \nu_{P_0}( \div(u) +W) \geq 0
\end{equation}
and
\begin{equation}  \label{NiXi14}
     \nu_{P_0}( w_m \tau_0) =  \nu_{P_0}( \div(w_m) +W) =m   \quad \for \quad m=0,1,... \; .
\end{equation}
We consider the sequence $(\dot{v}_j)_{j \geq 0} $  \eqref{NiXi11a}.
By \eqref{Th4-50}, $(\dot{v}_j)_{j = 0}^{m-1} $ is an $\FF_b$ linear basis of $L_m$.
 Let
\begin{equation}\label{T3-13b}
  V_j = \{ \dot{v}_j +\sum_{k=0}^{j-1} b_k \dot{v}_k \; | \; b_k \in \FF_b,\; k\in [0,j)  \}, \;\; \alpha(j) = \max_{v \in V_j} \;\nu_{P_0}(v \tau_0).
\end{equation}
It is easy to verify that $\alpha(j) \neq \alpha(j)$ for $i \neq j$.
We construct a sequence $(\ddot{v}_j)_{j \geq 0} $    as follows :
\begin{equation}\label{T3-14}
 \ddot{v}_0=\dot{v}_0, \quad \ddot{v}_j \in \{v\in V_j \;| \; \nu_{P_0}(v\tau_0) =\alpha(j)  \}, \quad \;\; j=1,2,... \;.
\end{equation}
It is easy to see that  $(\ddot{v}_j)_{j \geq 0} $ satisfy  the condition   \eqref{Th4-50}.
Bearing in mind (\ref{NiXi13})-(\ref{T3-13b}) and that $\ddot{v}_j \in L_m$ for $j <m$,  we get
\begin{equation}\label{T3-15}
  \nu_{P_0}(\ddot{v}_j\tau_0)  \neq \nu_{P_0}(\ddot{v}_k \tau_0)    \; \for \; j \neq k,
  \; \ad \; \nu_{P_0}(\ddot{v}_j \tau_0) =\alpha(j) \geq 0, \;\; j\geq 0.
\end{equation}
Hence, for all $f \in L_m $, we have
\begin{equation} \nonumber %   \label{T3-16}
  \nu_{P_0}(f \tau_0 ) \in \{\alpha(0),\alpha(1),... \} =: \dot{H}.
\end{equation}
%
%Now we will show what $\dot{H} =\NN_0 :=  \{n \; | \; n \geq 0\} $.
%By \eqref{NiXi10} and \eqref{NiXi11}, we have that  $w_j \in L_{j+g+1}$ with $\nu_{P_0}
Taking into account  \eqref{NiXi14} and \eqref{T3-15}, we obtain
\begin{equation}\label{T3-17}
 \dot{H} = \{n \; | \; n \geq 0\}= \NN_0.
\end{equation}
Suppose that $\alpha(j) > j+g$. By  \eqref{Th4-50} - \eqref{NiXi10}, $\ddot{v}_j \in  L_{j+1} =\cL((j-g +2)P_{s+1}+W) $. Hence $\ddot{v}_j \in \cL(X) $, with $X=(j-g+2)P_{s+1}+W -(j+g+1)P_0$.\\
Bearing in mind that $\deg(P_0) = \deg(P_{s+1})=1$ and $\deg(W) =2g-2$, we get $\deg(X) =-1$. Therefore  $\cL(X) =\{0\}$ and we have a contradiction.
%Bearing in mind \eqref{NiXi11a}..., we obtain
Hence
\begin{equation}\label{T3-13a}
 \alpha(j)  \leq j+g.
 \end{equation}
By \eqref{T3-17}, we have that  for every integer $k \geq 0$ there exists $r \geq 0$ with $\alpha(r)=k$. Therefore the map $\alpha: \; \NN_0 \to \NN_0 $ is an isomorphism.
  Hence there exist integers $\beta(k) \geq 0$ such that
% and  $\gamma_k\geq 0$ $(k=0,1,...)$
\begin{equation}  \label{PrT4-0}
   \beta(k)=\alpha^{-1}(k), \;    \alpha (\beta(k))  = k \;\; \ad \;\; \beta(\alpha(k))  = k  \;\for \; k=0,1,... \; .
\end{equation}
From \eqref{T3-13a}, we get for $j = \beta(k)$
\begin{equation} \label{PrT4-2}
     k = \alpha (\beta(k)) = \alpha (j) \leq j+g = \beta(k) +g.
\end{equation}
Let
\begin{equation} \label{PrT4-4}
 B_j =\{ r \geq 0  \; | \;  \alpha(r) < j  \}.
\end{equation}
Taking $r=\beta(k)$, we get $\alpha(r) =k$ and
\begin{equation} \label{PrT4-3}
          B_{j} = \{\beta(0),\beta(1),...,\beta(j-1) \}  \quad  \for \quad j \geq 1.
\end{equation}
Suppose  $  j \notin B_{j+g+1}$ for some $j$, then $j=\beta(j +g+l) $
for some  $l \geq 1$.
Using \eqref{PrT4-2}  with $k=j+g+l$, we obtain
\begin{equation} \nonumber
j+l =  (j +g+l)-g \leq \beta(j +g +l)=j.
\end{equation}
We have a contradiction. Hence
\begin{equation} \nonumber %\label{PrT4-3a}
 j \in  B_{j+g+1} \quad  \for \; {\rm all} \; \;j \geq 0.
\end{equation}
We consider the local expansion \eqref{NiXi05}, applied to $i=0$ :
\begin{equation} \label{NiXi05a}
  \dot{v}_r   \tau_0\; := \; \sum_{j = 1}^{\infty}
	\dot{c}^{(0)}_{j,r} t_0^{j-1}  , \;\; \where \;\;
	\dot{c}^{(0)}_{j,r} \in \FF_b, \;  j \geq 1, \;\;  \dot{C}^{(0)} =(\dot{c}^{(0)}_{j,r}  )_{j-1,r \geq 0}.
\end{equation}
%
%Let $\tilde{C}^{(0)} =(\dot{c}^{(0)}_{j,r}  )_{j-1,r \geq 0}$
Let $(x_n^{(0)}(\dot{C}^{(0)}))_{n \geq 0}$ be the digital sequence generated by the matrix $\dot{C}^{(0)}$.

Now we consider the matrix $\tilde{C}^{(i)} =(\dot{c}^{(i)}_{j,r}  )_{j-1,r \geq 0}$, obtained from
equation \eqref{NiXi05} and \eqref{NiXi05a}, where we take $\ddot{v}_r$ instead of  $\dot{v}_r$ $(i=0,1,...,s)$.
 Using  Theorem H, we obtain that $(x_n^{(0)}(\tilde{C}^{(0)}),\bx_n(\tilde{C}))_{n \geq 0}$ is the digital  $(t, s+1)$-sequence  with $t=g$.
Therefore we have proved the following lemma : \\  \\
{\bf Lemma 2.} { \it There exists a sequence  $(\ddot{v}_{j})_{j \geq 0}$ such that
$(x_n^{(0)}(\tilde{C}^{(0)}),\bx_n(\tilde{C}))_{n \geq 0}$ is the digital  $(t, s+1)$-sequence  with $t=g$ and $ \{0,1,...,m-1\} \subset B_{m+g} $.}\\ \\
In \S4.4, we will prove\\ \\
{\bf Theorem 4.} {\it  With the notations as above,
the set $[0,\gamma_1) \times ...\times [0,\gamma_s)$ is of bounded remainder  with respect to  $(\bx_n(\tilde{C}))_{n \geq 0}$ if and only if (\ref{Cond}) is true.} \\

\section{Proof} Consider the following condition
\begin{equation} \label{Prof0}
  \liminf_{m \to \infty} ( m - \bT(m))=\infty.
\end{equation}
We will prove \eqref{Prof0} for the generalized Halton sequence in \S4.3. For other considered sequences,
 assertion \eqref{Prof0} follows from Theorem D, Theorem E and Theorem H.

The sufficient part  of all considered theorems follows from Definition 2 and \eqref{Prof0}.
% (about properties of Halton : weak admissible ...).
Therefore we need only consider the case of necessity.  \\  \\
{\bf 4.1 Generalized Niederreiter sequence.  Proof of Theorem 1.} \\
From Theorem D, we have that $(\bx_n)_{n \geq 0}$ is the uniformly distributed digital
 weakly admissible $(t,s)$-sequence in base $b$. By Theorem B, in order to prove Theorem 1, we need only to check  condition (\ref{End}). By ( Le2, p.26, ref 4.6 ), we get
 \begin{equation}\label{gNied}
  y_{n,j}^{(i)}  =
     \underset{P_{\infty},x^{-1}}\Res \Big(\frac{y_{i,l,k(i,j)}(x)}
{p_i(x)^l}   \sum_{r=0}^{m-1}
		\bar{a}_{r}(n) x^{r+2}  \Big) = \underset{P_{\infty},x^{-1}}\Res \Big(\frac{y_{i,l,k(i,j)}(x)}
{p_i(x)^l}  n(x)  \Big)
\end{equation}
 \begin{equation} \nonumber
\with \; l =Q(i,j)+1 ,\; n(x) = \sum_{j=0}^{m-1} \bar{a}_{j}(n) x^{j+2} \; \ad \; \bar{a}_{j}(n) =  \underset{P_{\infty},x^{-1}}\Res (n(x) x^{-j-1} )  .
\end{equation}

% $ \quad	\quad {\rm  for \; all } \;\; j \in [1,d_i], \; i \in [1,s] .$
%
%
We take   $\dot{y}_{i,j,k}(x) =x^m y_{i,j,k}(x)$ instead of $y_{i,j,k}(x) $.
Now using Theorem D, we obtain from (\ref{gNied}), (\ref{Ap301}) - (\ref{Ap303}) that $(\dot{\bx}_n)_{0 \leq n < b^{\dot{m}}}$ is a $(t,\dot{m},s)$ net  for  $\dot{m} =s\tau_m + t$ with $\dot{x}_{n,j}^{(i)} =\phi^{-1}(\dot{y}_{n,j}^{(i)}) = x_{b^m n,j}^{(i)}$.
Bearing in mind that  $\dot{\bx}_{n}  = \bx_{b^m n}$,  we obtain (\ref{End}).
Hence  Theorem 1 is proved.  \qed  \\

{\bf 4.2 Xing-Niederreiter sequence.   Proof of Theorem 2.} \\
By Theorem B and Theorem E, in order to prove Theorem 2, we need only to check  condition (\ref{End}).

From \eqref{Ap300} - \eqref{Ap303}, we get that in order to obtain \eqref{End}, it suffices to prove that
\begin{align}
  &\# \{  n \in [0,b^M) \; |  \; y^{(i)}_{n,j} = u^{(i)}_{j},\;  j \in [1,  \tau_m] \;\;\for\; i \in [1,s],\;
    \ad \; a_{j-1}(n)= u^{(0)}_{j}    \nonumber \\
  &    j \in [1,m]    \} > 0, \qquad \;\; \with \quad M=s\tau_m +(m+2g)(2g+1)e_0 +m_0, \label{T2Pr-1}
\end{align}
$m_0 =2g+2   + e_1+\cdots +e_s $, for all $u^{(i)}_{j} \in \FF_b$. \\
Let
\begin{equation} \nonumber
   \delta(\fT) =   \begin{cases}
    1,  & \; {\rm if}  \;  \fT  \;{\rm is \;true},\\
    0, &{\rm otherwise}.
  \end{cases}
\end{equation}
Let $k^{(0)}_{j+1} =z_{h(j)} = z^{h(j)}$ for $j \in H_1$ with
 $H_1 =\NN_0 \setminus H_2 =\{h(0),h(1),...\}$, $H_2 = \{n_0,n_1,...,n_g   \}$.
 %and let    $H_3 =h_1(H_1)$, $H_4 =h_1(H_2)$.
 From  \eqref{XiNi10}, we have
 \begin{equation} \nonumber % \label{T2Pr-1a}
   a^{(0)}_{j,r} =\delta(j-1 =r \in H_1), \quad j \geq 1.
\end{equation}
Let $ c^{(0)}_{j,r} :=a^{(0)}_{j,h(r)} $. By  \eqref{Ap301}, \eqref{Ap302} and \eqref{XiNi07}, we get
%\eqref{XiNi08}, \eqref{XiNi10} and \eqref{XiNi07}, we get
%
\begin{align}
&   c^{(0)}_{j,r} = \delta(j-1 =h(r)), \; \quad   \; y^{(0)}_{n,j} = \sum_{r \geq 0} \bar{a}_{r}(n)c^{(0)}_{j,r}= \bar{a}_{h(j-1)}(n), \;\;    \label{T2Pr-2}      \\
&   x^{(0)}_n =\sum_{j \geq 1} a_{h(j-1)}(n)/b^j \quad \ad \;\; k^{(0)}_j =z^{h(j-1)} \in  \cL\big( h(j-1)(2g+1) P_0 \big), \; j \geq 1. \nonumber
\end{align}
So, we obtain a digital $s+1$-dimensional sequence   $(x^{(0)}_n ,\bx_n)_{n \geq 0}$.

Let $n =\sum_{r=0}^{M-1} a_r(n) b^r$ and let
\begin{equation}   \nonumber %\label{T2Pr-3}
\dot{n} =\sum_{r \in H_1} a_{r}(n) b^{r}, \; \ddot{n} =\sum_{r \in H_2} a_{r}(n) b^{r}, \;
 \dot{U}= \{  \dot{n}  | n \in [0,b^M) \}, \\
  \ddot{U}= \{  \ddot{n}  | n \in [0,b^M) \}.
\end{equation}
By \eqref{Ap301}, \eqref{XiNi12} and \eqref{T2Pr-2}, we get
\begin{align} \nonumber %  \label{T2Pr-4}
& y^{(i)}_{n,j} = \sum_{  r \geq 0} \bar{a}_{r}(n) c^{(i)}_{j,r}=
\sum_{  r \in H_1} \bar{a}_{r}(n) c^{(i)}_{j,r}+ \sum_{  r \in H_2} \bar{a}_{r}(n) c^{(i)}_{j,r}
= y^{(i)}_{\dot{n},j} + y^{(i)}_{\ddot{n},j} ,  \nonumber \\
& \quad  i \in [1,s], \quad  y^{(0)}_{n,j} = y^{(0)}_{\dot{n},j} =\bar{a}_{h(j-1)}(n) ,\;\; y^{(0)}_{\ddot{n},j} =0, \quad j \geq 1.  \nonumber
\end{align}
We fix $\tilde{n} \in \ddot{U}$. Let
\begin{align} \nonumber %  \label{T2Pr-6}
& A_{\bu,\tilde{n} } = \{  n \in [0,b^M) \; | \;  y^{(i)}_{n,j} = u^{(i)}_{j},\;  j \in [1,  \tau_m], \; i \in [1,s], \nonumber \\
&    \;\;  y^{(0)}_{n,j} = u^{(0)}_{j},\; j\in [1,m],\quad  \ddot{n} = \tilde{n}    \}.
\end{align}
It is easy to verify that statement \eqref{T2Pr-1}  follows from the next assertion
\begin{equation} \label{T2Pr-7}
     \#A_{\bu,\tilde{n} } >0 \quad  \forall \; u^{(i)}_{j} \in \FF_b, \quad
     \tilde{n} \in \ddot{U}.
\end{equation}
Taking into account that $y^{(i)}_{n,j} = y^{(i)}_{\dot{n},j} + y^{(i)}_{\ddot{n},j}$, we get
\begin{equation} \nonumber  %\label{T2Pr-8}
A_{\bu,\tilde{n} } = \{ \dot{n} \in \dot{U} \; | \;  y^{(i)}_{\dot{n},j} = \dot{u}^{(i)}_{j},\;  j \in [1,  \tau_m], \; i \in [1,s], \; y^{(0)}_{\dot{n},j} = u^{(0)}_{j},\; j\in [1,m]   \},
\end{equation}
where $ \dot{u}^{(i)}_{j} = \dot{u}^{(i)}_{j} - y^{(i)}_{\tilde{n},j}  $.
%
%According to Theorem F, we have that $(\bx_n)_{n \geq 0}$ is a  uniformly distributed digital  weakly admissible $(t,s)$ sequence in base $

According to \eqref{Ap301}, \eqref{XiNi12} and \eqref{T2Pr-2},
in order to prove \eqref{T2Pr-7} , it suffices to show that the vectors
\begin{equation} \label{T2Pr-8}
\pi_M(\bc^{(i)}_j) = (c^{(i)}_{j,0}, ...,c^{(i)}_{j,M-1} ) \in F^M_b,
    \quad \with
 \quad  1 \leq j \leq d_i, \; 0 \leq i \leq s,
\end{equation}
$d_i =\tau_m $,  $ 1 \leq i \leq s$ and  $d_0 =m $, are linearly independent over $\FF_b$.
% where $M=s\tau_m +(m+2g)(2g+1)e_0 +5g$.

To prove this statement, we closely follow [DiPi, p.282].
 Suppose that we have
\begin{equation} \nonumber %   \label{T2Pr-9}
   \sum_{j=1}^{ m} f^{(0)}_j  \pi_M(\bc^{(0)}_j) +\sum_{i=1}^s \sum_{j=1}^{\tau_m} f^{(i)}_j  \pi_M(\bc^{(i)}_j)  = \bs \in F^M_b
\end{equation}
for some $f^{(i)}_j  \in \FF_b$ with $ \sum_{j=1}^{ m} |\phi^{-1}(f^{(0)}_j)|+\sum_{i=1}^s \sum_{j=1}^{\tau_m} |\phi^{-1}(f^{(i)}_j)|   >0$.\\
 We put $f^{(0)}_{r}=0$ for $r>m$. Hence
\begin{equation} \nonumber %     \label{T2Pr-10}
   \sum_{j=1}^{ m} f^{(0)}_j  c^{(0)}_{j,r} +\sum_{i=1}^s \sum_{j=1}^{\tau_m} f^{(i)}_j  c^{(i)}_{j,r}=0  \quad \for \quad r \in [0,M).
\end{equation}
By   \eqref{XiNi12} and \eqref{T2Pr-2}, we obtain $ c^{(i)}_{j,r} =a^{(i)}_{j,h(r)} $  for $ 1 \leq i \leq s$ and $ c^{(0)}_{j,r}=\delta(j-1\\=h(r))$. Therefore
\begin{equation} \label{T2Pr-11}
 0=\sum_{j=1}^{ m} f^{(0)}_j  \delta(j-1=h(r)) +\sum_{i=1}^s \sum_{j=1}^{\tau_m} f^{(i)}_j  a^{(i)}_{j,h(r)}= f^{(0)}_{h(r)+1} + \sum_{i=1}^s \sum_{j=1}^{\tau_m} f^{(i)}_j  a^{(i)}_{j,h(r)}
\end{equation}
$\for \; r \in [0,M)$. \\
 Now consider the element $ \alpha \in \FF_b$  given by  $\alpha= \alpha_1 + \alpha_2$, where
\begin{equation} \label{T2Pr-12}
\alpha_1=\sum_{r=0}^{m-1} f^{(0)}_{h(r)+1} z_{h(r)} , \quad \alpha_2 =\sum_{i=1}^s \sum_{j=1}^{\tau_m} f^{(i)}_j  k^{(i)}_j -  \sum_{i=1}^s \sum_{j=1}^{\tau_m} f^{(i)}_j   \sum_{u=0}^g a^{(i)}_{j,n_u}  w_u.
\end{equation}
%Put $H_2=\{ n_0,n_1,...,n_g    \}$  and
Using \eqref{XiNi10}, we get
\begin{equation} \nonumber
 \alpha_2 =  \sum_{i=1}^s \sum_{j=1}^{\tau_m} f^{(i)}_j \Big( \sum_{r=0}^{\infty}  a^{(i)}_{j,r} z_r -  \sum_{u=0}^g a^{(i)}_{j,n_u}  z_{n_u} \Big) =
   \sum_{r \in H_1} \Big(  \sum_{i=1}^s \sum_{j=1}^{\tau_m} f^{(i)}_j    a^{(i)}_{j,r} \Big) z_r .
\end{equation}
From \eqref{XiNi12}, (\ref{T2Pr-11}) and \eqref{T2Pr-12}, we obtain
\begin{equation} \nonumber
 \alpha=  \sum_{r \geq 0} \Big( f^{(0)}_{h(r)+1}+ \sum_{i=1}^s \sum_{j=1}^{\tau_m} f^{(i)}_j    a^{(i)}_{j,{h(r)}} \Big) z_{h(r)} = \sum_{r \geq M} \Big( f^{(0)}_{h(r)+1}+ \sum_{i=1}^s \sum_{j=1}^{\tau_m} f^{(i)}_j    a^{(i)}_{j,{h(r)}} \Big) z_{h(r)} .
\end{equation}
%
%\begin{equation} \nonumber\alpha_1= \sum_{r=0}^{m-1} f^{(0)}_{r+1}     z_{r},\;\; \ad \;\;  f^{(0)}_{r+1} + \sum_{i=1}^s \sum_{j=1}^{\tau_m} f^{(i)}_j    a^{(i)}_{j,r} =0 \;\; \for \;\; r\in [0,M+g] \setminus H_2\end{equation}
%and $f^{(0)}_{r+1}=0$ for $r \in R$.
%for the first $M$ non-negative integers $r$ that are not in $R$.
Hence
\begin{equation} \label{T2Pr-16}
       \nu_{P_{\infty}}(\alpha)  \geq M.
\end{equation}
Furthermore, \eqref{XiNi06}, \eqref{XiNi08}, \eqref{XiNi07}, \eqref{T2Pr-2}  and \eqref{T2Pr-12} yield
\begin{equation} \label{T2Pr-17}
    \alpha_1 \in  \cL( (m+2g)(2g+1) P_0), \;\;   \alpha_2 \in  \cL\Big( G +\sum_{i=1}^s([\tau_m/e_i]+1) P_i  \Big).
\end{equation}
Combining \eqref{T2Pr-16} and \eqref{T2Pr-17}, we obtain
\begin{equation} \nonumber
   \alpha \in    \cL\Big( G +   \sum_{i=1}^s ([\tau_m/e_i]+1) P_i  +(m+2g)(2g+1)P_0  - M P_{\infty}      \Big).
\end{equation}
But  from \eqref{T2Pr-1}, we have
\begin{align} \nonumber
 &      \deg\Big(  G +   \sum_{i=1}^s ([\tau_m/e_i]+1) P_i  +(m+2g)(2g+1)P_0  - M P_{\infty} \Big)  \nonumber \\
 &     = 2g +     \sum_{i=1}^s ([\tau_m/e_i]+1)e_i  +(m+2g)(2g+1) e_0  - M  \nonumber      \\
 & \leq 2g +s\tau_m + e_1+\cdots +e_s +(m+2g)(2g+1) e_0 -M <0.     \nonumber
\end{align}
Hence
\begin{equation} \nonumber
       \cL\Big(G +   \sum_{i=1}^s ([\tau_m/e_i]+1) P_i  +(m+2g)(2g+1)P_0  - M P_{\infty}    \Big) =\{ 0 \}
\end{equation}
by \eqref{No20} and therefore we have $\alpha=0$.

By \eqref{XiNi06}, we have $\nu_{P_0}(k^{(i)}_j) \geq 0$ and $\nu_{P_0}(w_u) \geq 0$ for all $i,j,u$.
According to \eqref{T2Pr-12}, we get $\nu_{P_0}(\alpha_2) \geq 0.$
Suppose that $\alpha_1 \neq 0$.
Taking into account that $z_0 =z_{n_0} =w_0 \neq z_{h(r)} $ for $r\geq 0$, we obtain from \eqref{T2Pr-12} that $\nu_{P_0}(\alpha_1) < 0$. We have a contradiction. Hence $\alpha_1 = 0$ and $\alpha_2 = 0$. From Lemma C, we conclude that  $f^{(i)}_j =0$ for all
$i,j$.
% From \eqref{T2-3} and Lemma C we conclude that  $ f^{(i)}_j =0$  for $  1 \leq i \leq s, j \geq 1$.
Hence the system \eqref{T2Pr-8} is linearly independent over $\FF_b$.

Thus  \eqref{T2Pr-1} is true and $ (\bx_n)_{ n \geq 0}$ satisfies the condition (\ref{End}).
 By Theorem~E,  $ (\bx_n)_{ n \geq 0}$
 is the $d-$admissible uniformly distributed digital $(t,s)$-sequence in base $b$.
 Applying  Theorem B, we get the assertion of Theorem 2. \qed  \\

{\bf 4.3 Generalized Halton-type sequence.  Proof of Theorem 3.}\\ \\
{\bf Lemma 3.} { \it The sequence $(\bx_n)_{n \geq 0}$ is   uniformly distributed in $[0,1)^s$ .}\\ \\
{\bf Proof.} By Lemma A, in order to prove Lemma 3, it  suffices to show that $m -T(m) \to \infty$ for
$m  \to \infty$.
Let $R_k = \max_{1 \leq i \leq s} n_{i,k}$, $k=1,2,...$  . We define $j_{i,k}$ from the following condition
$ n_{i,j_{i,k}} \geq R_k > n_{i,j_{i,k}-1} $.
Let $ \tilde{R}_k =\sum_{i=1}^s n_{i,j_{i,k}}  $.

We consider the definition of $(t,m,s)$ net.
% with $d_1+d_2 +\cdots + d_s =R_k$, $m \geq \tilde{R}_k$ and
%$t=m-R_k$. It is easy to see that.
Suppose that for all
$$
 E = \prod_{i=1}^s [a_ib^{-d_i},(a_i+1)b^{-d_i}),\;\;
  \with \; a_i =\sum_{j=1}^{d_i} a_{i,j}b^{j-1}, \;a_{i,j} \in \ZZ_b, \;  d_i \geq 0,
$$
$ 1 \le i \leq s, \; d_1+ \cdots + d_s =R_k  $, we have
\begin{align}
& \# \{ n \in [0,b^m)\; |\; x_n \in E \} =\# \{  n \in [0,b^m) \; | \;  y^{(i)}_{n,j} = u^{(i)}_{j},\;  j \in [1,  d_i],  i \in [1,s]     \}   \nonumber\\
& =b^{m-R_k}, \quad \where \quad
   m \geq \tilde{R}_k +(3g+3)e_0, \; u^{(i)}_{j} = \phi^{-1}(a_{i,j})\in \FF_b,
      \label{T3Pr-10}
\end{align}
 $j \in [1,  d_i], \; i \in [1,s]$.

By Definition 2, we get that $(\bx_n)_{n \geq 0}$ is a $(T,s)$-sequence in base $b$ with $m-R(k) \geq T(m)$ for $ m  \geq \tilde{R}_k +(3g+3)e_0 $.
Bearing in mind that $R(k) \to \infty$ for $k \to \infty$, we obtain the assertion of Lemma 4.

Taking into account that  $d_i \leq R_k \leq n_{i,j_{i,k}}$ for $ 1 \leq i \leq s $, we get that in order to prove \eqref{T3Pr-10}, it suffices to verify that
\begin{equation} \label{T3Pr-11}
\# \{  n \in [0,b^m) \; | \;  y^{(i)}_{n,j} = u^{(i)}_{j},\;  j \in [1,  n_{i,j_{i,k}}], \; i \in [1,s]     \} =b^{m-\tilde{R}_k}
\end{equation}
for all $u^{(i)}_{j} \in \FF_b$, with $j \in [1, n_{i,j_{i,k}}], \; i \in [1,s]$.\\
%
%
% \tilde{R}_k}   \eqref{T3-11} and \eqref{T3-12}
%
Let $ \cM=(m_0e_0  +2g-1)P_{\infty}  $ with $m_0=[m/e_0] -2g-1$.\\
By Lemma 1, we obtain that there exist sets $H_1$ and $H_2$ such that $H_1 \cup H_2 =\{0,1,...,m-1\}$,  $H_1 \cap H_2 = \emptyset$,  $(w^{(0)}_r)_{r \in H_1}$  is the $\FF_b$ linear basis of $\cL(\cM)$ and $\#H_2 =m-m_0e_0-g=:g_1 $, with $g_1 -e_0(2g+1)-g \in [0,e_0)  $.
 Let $n =\sum_{r=0}^{m-1} a_r(n) b^r$ and let
\begin{equation} \nonumber %   \label{T3Pr-12}
\dot{n} =\sum_{r \in H_1} a_{r}(n) b^{r}, \; \ddot{n} =\sum_{r \in H_2} a_{r}(n) b^{r}, \;
 \dot{U}= \{  \dot{n}  | n \in [0,b^m) \}, \; \ddot{U}= \{  \ddot{n}  | n \in [0,b^m).
\end{equation}
So
\begin{equation} \nonumber %   \label{T3Pr-13}
 f_n = \sum_{  r = 0}^{m-1} \bar{a}_{r}(n) w^{(0)}_r \in \cL(\cM) \quad \Longleftrightarrow \quad
   n=\dot{n}, \;\; \for \; n \in [0,b^m) .
\end{equation}
We fix $\tilde{n} \in \ddot{U}$. Let
\begin{equation}  \nonumber % \label{T3Pr-14}
A_{\bu,\tilde{n} } = \{  n \in [0,b^m) \; | \;  y^{(i)}_{n,j} = u^{(i)}_{j},\;  j \in [1,  n_{i,j_{i,k}}], \; i \in [1,s], \;\;  \ddot{n} = \tilde{n}    \}.
\end{equation}
It is easy to see that statement \eqref{T3Pr-11}  follows from the next assertion
\begin{equation} \label{T3Pr-15}
     \#A_{\bu,\tilde{n} } = b^{m-\tilde{R}_k -g_1} \quad  \forall \; u^{(i)}_{j} \in \FF_b, \quad
     \tilde{n} \in \ddot{U}.
\end{equation}
Taking into account that $y^{(i)}_{n,j} = y^{(i)}_{\dot{n},j} + y^{(i)}_{\ddot{n},j}$, we get
\begin{equation} \nonumber % \label{T3Pr-16}
A_{\bu,\tilde{n} } = \{ \dot{n} \in \dot{U} \; | \;  y^{(i)}_{\dot{n},j} = \dot{u}^{(i)}_{j},\;  j \in [1,  n_{i,j_{i,k}}], \; i \in [1,s]   \},
\end{equation}
where $ \dot{u}^{(i)}_{j} = \dot{u}^{(i)}_{j} - y^{(i)}_{\tilde{n},j}  $.
Let
\begin{equation} \nonumber % \label{T3Pr-17}
   \hat{\psi}(f_n):=(    y^{(1)}_{n,1},...,  y^{(1)}_{n,n_{1,j_{1,k}}}  , ..., y^{(s)}_{n,1},...,  y^{(s)}_{n,n_{s,j_{s,k}}}           )  \in \FF_b^{\tilde{R}_k}.
\end{equation}
We consider the map  $\breve{\psi} \; : \;\cL(\cM) \to \FF_b^{\tilde{R}_k}$
defined by
\begin{equation} \nonumber %   \label{T3Pr-18}
   \breve{\psi}(\dot{f}): = \hat{\psi}(f_n)    \quad {\rm where} \quad  \cL(\cM) \ni \dot{f}=f_n \quad \with \; {\rm some} \; n  \in \dot{U}.
\end{equation}
Note that $\breve{\psi}$ is a linear transformation between vector spaces over $\FF_b$.
It is clear that in order to prove \eqref{T3Pr-15}, it suffices to verify that $\breve{\psi}$ is
surjective.
To prove this, it is enough to show that
\begin{equation} \label{T3Pr-19}
   \dim\big( \cL(\cM)/ \ker(\breve{\psi})   \big) =\tilde{R}_k.
\end{equation}
%
%Now we will prove that $\breve{\psi}$ is surjective.
%To prove this, it suffices to show that
%
Using \eqref{T3-1}, \eqref{T3-3} and \eqref{T3-12},   we get that $w^{(i)}_l \equiv 0 \; (\mod \; \cP_{i,j_{i,k}})$ for\\ $l \geq n_{i,j_{i,k}}$.
By \eqref{T3-1}, \eqref{T3-3},  \eqref{T3-12}, and  \eqref{T3Pr1}, we derive that  $ y^{(i)}_{n,j}=0$ for all $j \in [1, n_{i,j_{i,k}}]$  if and only if $\dot{f}=f_n  \equiv \; 0\; (\mod \; \cP_{i,j_{i,k}})$  for $i \in [1,s]$.\\
From the definition of $\breve{\psi}$ it is clear that
\begin{equation} \nonumber %  \label{T3Pr-4}
    \ker(\breve{\psi}) =  \cL(H), \quad \with \quad H= \cM -\sum_{i=1}^s  \cP_{i,n_{i,j_{i,k}}}.
\end{equation}
Using Riemann-Roch's theorem, we obtain that $\dim(\cM)=m_0e_0+g=m-g_1$, where
$\deg(\cM) =m_0e_0+2g-1$ and
\begin{equation} \nonumber %     \label{T3Pr-5}
   \deg(H)  =  m_0e_0+2g-1 -\sum_{i=1}^s n_{i,j_{i,k}} =m+g-1 -g_1 -\tilde{R}_k.
\end{equation}
Hence  $\dim(\ker(\psi))=m-\tilde{R}_k-g_1\geq (3g+3)e_0-g_1 \geq (3g+3)e_0 -(2g+2)e_0 -g      \geq 1$, $\dim(\cM)=m-g_1$  and   \eqref{T3Pr-19} follows.
So $\breve{\psi}$ is indeed surjective. Therefore \eqref{T3Pr-15} and Lemma 3 are proved. \qed \\ \\
{\bf Lemma 4.} { \it  The sequence $(\bx_n)_{n \geq 0}$  satisfies  condition \eqref{End}.}\\  \\
{\bf Proof.}  Let
\begin{equation}      \label{T3Pr-0}
M =([M_1/e_0]+3g+1)e_0 , \;\;  M_1=\sum_{i=0}^s n_{i,j_{i,m}} \; \where \; n_{j_{i,m}} \geq \tau_m > n_{j_{i,m}-1}
\end{equation}
for $i \in [1,s]$,  $n_{0,j_{0,m}} = ( [m/e_0] +1)e_0$  $j_{0,m} =[m/e_0] +1$. \\
Bearing in mind that $y^{(0)}_{n,j}=\bar{a}_{j-1}(n) $, $(j=1,2,...)$,
 we get from \eqref{T3Pr2} - \eqref{T3Pr4}, that in order to obtain \eqref{End}, it  suffices to prove that
\begin{equation}  \label{T3Pr-1}
  \# \{  n \in [0,b^M) \; | \;  y^{(i)}_{n,j} = u^{(i)}_{j},\;\;\;\;  j \in [1,   n_{i,j_{i,m}}] \;\;\;\; \for \;\;  i \in [0,s]      \} > 0
\end{equation}
for all $u^{(i)}_{j} \in \FF_b$.
%
% By \eqref{T3-9}-\eqref{T3-11}  we get that
%
Let $ \cM=(([M_1/e_0]+1)e_0+2g-1)P_{\infty}$. By \eqref{T3-0}, $\deg(P_{\infty})=1$. Hence $\deg(\cM)=([M_1/e_0] +1)e_0+2g-1 $. Using Riemann-Roch's theorem, we obtain that
\begin{equation} \label{T3Pr-2}
    \dim(\cM)=([M_1/e_0]+1)e_0+g =M_1+g_1+g \quad \with \quad  g_1:= ([M_1/e_0]+1)e_0 - M_1    .
\end{equation}

By Lemma 1, we get that an $\FF_b$ linear basis of $\cL(\cM)$ can be chosen from the set $\{ w^{(0)}_{0},..., w^{(0)}_{M-1}\}$ with $M = ([M_1/e_0] +3g+1)e_0 =n_{0,[M_1/e_0] +3g+1}$.

 Let $n =\sum_{r=0}^{M-1} a_r(n) b^r$ and let $ f_n = \sum_{r=0}^{M-1}  \bar{a}_{r}(n) w^{(0)}_r$.
We get that for all $ \dot{f}  \in \cL(\cM)$ there exists  $ n \in [0,b^M)$ such that $\dot{f}=f_n$.

From \eqref{T3Pr1},  we have
\begin{equation} \label{T3Pr-3a}
   f_n = \sum_{  j =1}^{\infty}   y^{(i)}_{n,j} w^{(i)}_{j-1}, \qquad  0 \leq i \leq s.
\end{equation}
Let
\begin{equation} \label{T3Pr-3}
   \psi(f_n):=(    y^{(0)}_{n,1},...,  y^{(0)}_{n,n_{0,j_{0,m}}}, ..., y^{(s)}_{n,1},...,  y^{(s)}_{n,n_{s,j_{s,m}}}           )  \in \FF_b^{M_1}.
\end{equation}
Consider the map  $\dot{\psi} \; : \;\cL(\cM) \to \FF_b^{M_1}$
defined by
\begin{equation} \nonumber % \label{T3Pr-3b}
   \dot{\psi}(\dot{f}): = \psi(f_n)    \quad {\rm where} \quad  \dot{f}=f_n \quad \with \; {\rm some} \; n \in [0,b^{M}).
\end{equation}
We see that in order to obtain \eqref{T3Pr-1}, it  suffices to verify that $\dot{\psi}$ is
surjective. \\
To prove this, it suffices to show that
\begin{equation} \label{T3Pr-3c}
   \dim\big( \cL(\cM)/ \ker(\dot{\psi})   \big) =M_1.
\end{equation}

Using \eqref{T3-1}, \eqref{T3-3} and \eqref{T3-12},   we get that $w^{(i)}_k \equiv 0 \; (\mod \;  \cP_{i,j_{i,m}})$ for $k > n_{i,m}$.
From \eqref{T3Pr-3a}, \eqref{T3-1}, \eqref{T3-3} and \eqref{T3-12}, we derive that  $ y^{(i)}_{n,j}=0$ for all $j \in [1, n_{i,j_{i,m}}]$  if and
only if $ f_n \equiv 0 \; (\mod \; \cP_{i,j_{i,m}})$  for $i \in [0,s]$. \\
From the definition of $\dot{\psi}$ it is clear that
\begin{equation} \nonumber %   \label{T3Pr-4a}
    \ker(\dot{\psi}) =  \cL(H), \quad \with \quad H= \cM -\sum_{i=0}^s  \cP_{i,j_{i,m}} .
\end{equation}
Using    \eqref{T3Pr-0}, \eqref{T3Pr-2}, \eqref{T3-0} and Riemann-Roch's theorem, we obtain that
\begin{equation}  \nonumber % \label{T3Pr-5a}
   \deg(H)  =  M_1 +g_1+2g-1 -\sum_{i=0}^s n_{i,j_{i,m}}  =g_1 +2g-1
\end{equation}
and $\dim(\ker(\psi))=g_1+g$. By \eqref{T3Pr-2}, $\dim(\cM)=M_1+g_1 +g$. Hence
$  \dim\big( \cL(\cM)/ \ker(\dot{\psi})   \big) =M_1$. Therefore \eqref{T3Pr-3c} is true.
So $\psi$ is indeed surjective and \eqref{T3Pr-1} follows. Therefore Lemma 4 is proved. \qed \\ \\
{\bf Lemma 5.} { \it The sequence $(\bx_n)_{n \geq 0}$ is   weakly admissible.} \\ \\
{\bf Proof.} Suppose that $x^{(i)}_{n} =x^{(i)}_{k}$ for some $i,n,k$. From \eqref{T3Pr-3} and
\eqref{T3Pr2}-\eqref{T3Pr4}, we get that
$y^{(i)}_{n,j} =y^{(i)}_{k,j}$ for $j \geq 1$.

Using \eqref{T3Pr-3a}, we have
\begin{equation} \nonumber
 f_n = \sum_{  j \geq 1} y^{(i)}_{n,j} (n) w^{(i)}_{j-1} .
\end{equation}
Hence $f_n =f_k$. Taking into account that $(w^{(0)}_{r})_{r \geq 0}$ is an $\FF_b$ linear basis of $O_{F}$, we obtain from \eqref{T3Pr0a}, that $n=k$.  By Definition 7, Lemma 5 is proved. \qed

Applying Theorem B, we get the assertion of Theorem 3. \qed \\ \\
{\bf 4.4 Niederreiter-Xing sequence. Proof of Theorem 4.} \\
Similarly to the proof of Lemma 5, we get that $(\bx_n)_{n \geq 0}$ is   weakly admissible.
By Lemma 2,  $(\bx_n)_{n \geq 0}$ is  the digital uniformly distributed sequence.

According to \eqref{Ap301}, \eqref{Ap302}, \eqref{End} and Theorem B, in order to prove Theorem 4,
 it is enough to verify that
\begin{equation} \label{T4}
  \# \{  n \in [0,b^M)  |   y^{(i)}_{n,j} = u^{(i)}_{j},\;  j \in [1,  \tau_m],   i \in [1,s] ,\;
 a_j(n)=0   \;\for\;j \in [0,m)   \} > 0
\end{equation}
for all $u^{(i)}_{j} \in \FF_b$, where  $M =s\tau_m +m +2g+2$.

Bearing in mind that by Lemma 2 $(x^{(0)}_n,\bx_n)_{n \geq 0}$ is a $(g,s+1)$ sequence, we obtain
\begin{equation} \nonumber
  \# \{  n \in [0,b^M)  |   y^{(i)}_{n,j} = u^{(i)}_{j},  j \in [1,  \tau_m],i \in [1,s],\;  y^{(0)}_{n,j}=0,
      j \in [1,m+g+2]     \} > 0
\end{equation}
for all $u^{(i)}_{j} \in \FF_b$.

Therefore, in order to prove \eqref{T4},  it  suffices to verify that
\begin{equation} \label{T4-1}
  {\rm if} \;  y^{(0)}_{n,j}=0
   \;\for\;    j \in [1,m+g+2]    \quad \then \quad   a_j(n)=0   \;\for \;j \in [0,m).
\end{equation}
Now we will prove \eqref{T4-1} :

  From \eqref{NiXi05} and \eqref{T3-14}, we have
 $\ddot{v}_r \dot{\tau}_0 =\sum_{j \geq 1} \dot{c}^{(0)}_{j,r} t_0^{j-1}$  with
$\nu_{P_0} (\ddot{v}_r \tau_0)\\  = \alpha(r)$. Hence
  $\dot{c}^{(0)}_{j,r} =0$ for $j \leq \alpha(r)$ and  $\dot{c}^{(0)}_{j,r} \neq 0$ for $j =\alpha(r)+1$.

Using \eqref{Ap301}, \eqref{PrT4-0} and  \eqref{PrT4-4} %and \eqref{PrT4-6},
 we obtain $\dot{c}^{(0)}_{j,\beta(j-1)} \neq 0$ and
\begin{equation} \nonumber  %\label{PrT4-1}
    y^{(0)}_{n,j}  =  \sum_{r \geq 0} \bar{a}_r(n) \dot{c}^{(0)}_{j,r}  =  \sum_{\alpha(r) <j} \bar{a}_r(n) \dot{c}^{(0)}_{j,r} = \sum_{r \in B_j} \bar{a}_r(n) \dot{c}^{(0)}_{j,r}, \;\; \quad j \geq 1.
\end{equation}
We apply  induction and consider the case $j=1$.
By \eqref{PrT4-3}, we see that $\bar{a}_{\beta(0)}(n) =0 $ if $ y^{(0)}_{n,1} =0 $.  Suppose that
 $\bar{a}_{\beta(0)}(n) =\cdots = \bar{a}_{\beta(l-1)}(n)= 0$
  if $ y^{(0)}_{n,1} = \cdots =  y^{(0)}_{n,l} =0 $ for some  $l \geq 1 $ .%\in [0, m+20g-1]$.
Now let $ y^{(0)}_{n,1} = \cdots =  y^{(0)}_{n,l} =y^{(0)}_{n,l+1} =0$.
We see
\begin{equation} \nonumber %\label{PrT4-2a}
   0=y^{(0)}_{n,l+1}   = \sum_{r \in B_{l+1}} \bar{a}_r(n) \dot{c}^{(0)}_{l+1,r} =
   \sum_{r \in B_{l+1} \setminus B_l} \bar{a}_r(n) \dot{c}^{(0)}_{l+1,r}= \bar{a}_{\beta(l)}(n) \dot{c}^{(0)}_{l+1,\beta(l)} .
\end{equation}
Bearing in mind that $\dot{c}^{(0)}_{l+1,\beta(l)} \neq 0$, we get $\bar{a}_{\beta(l)}(n)= 0$. \\
 Therefore if $y^{(0)}_{n,j} =0  $
 for all  $  1 \leq j  \leq m+g+1$, then $a_{\beta(j-1)}(n) =0$   for all  $  1 \leq j  \leq m+g+1$.
% \eqref{PrT4-3a}
Using Lemma 2, we get
 $a_r(n) =0$   for all  $  0 \leq r  \leq m-1$.

Hence \eqref{T4-1} is true  and Theorem 4 follows. \qed \\ \\
{\bf Aknowledgment.} Parts of this work were started at the Workshop "Discrepancy
Theory and Quasi-Monte Carlo methods" held at the Erwin Schr\"odinger Institute,
September 25 - 29, 2017. \\

{\bf Address}: Department of Mathematics,
Bar-Ilan University, Ramat-Gan, 5290002, Israel \\
{\bf E-mail}: mlevin@math.biu.ac.il\\

\end{document}